\documentclass{article}
\usepackage[utf8]{inputenc}
\usepackage{amsmath}
\usepackage{amsthm}
\usepackage{comment}
\usepackage{amsfonts} 
\usepackage{graphicx}
\usepackage{lscape}
\usepackage{tablefootnote}
\usepackage{algorithm}
\usepackage{algpseudocode}
\usepackage{array, booktabs}
\usepackage{dsfont} 
\usepackage[titletoc, title]{appendix}
\usepackage[round]{natbib}
\usepackage{color}
\numberwithin{equation}{section}
\usepackage{soul}
\usepackage{subfig}

\usepackage{arydshln}

\usepackage{amssymb}
\def\rev#1{{\color{black}{#1}}}

\date{\today}

\author{David Pinzon Ulloa\footnote{Corresponding author. \texttt{david.pinzon@umontreal.ca} \\ \texttt{Université de Montréal, Pavillon André Aisenstadt, 2920 chemin de la Tour, Montréal, QC H3T1J4 Canada.
 }} \footnote{\texttt{Escuela Superior Politécnica del Litoral, Campus Gustavo Galindo - km. 30.5 vía perimetral, Guayaquil - Ecuador}} \and Ammar Metnani\footnote{\texttt{ammar.metnani@umontreal.ca} \\\texttt{Université de Montréal, C.P. 6128, Succursale
Centre-ville, Montreal, QC H3C 3J7, Canada}} \footnote{\texttt{CIRRELT, Pavillon André Aisenstadt, bureau 3520
2920, chemin de la Tour
Université de Montréal
Montréa}} \and  Emma Frejinger\footnote{\texttt{emma.frejinger@umontreal.ca} \\\texttt{DIRO (FAS), Université de Montréal, C.P. 6128, Succursale
Centre-ville, Montreal, QC H3C 3J7, Canada}}}

\title{Usage of the \texttt{\textbackslash author} command}

\title{A Capacitated Collection-and-Delivery-Point Location Problem with Random Utility Maximizing Customers}

\begin{document}

\maketitle

\begin{abstract}
We consider a strategic decision-making problem where a logistics provider (LP) seeks to locate collection and delivery points (CDPs) with the objective \rev{of} reduc\rev{ing} total logistics costs. The customers maximize utility \rev{based on} their perception of home delivery service as well as the characteristics of the CDPs, including their location. At the strategic planning level, the LP does not have complete information about customers' preferences and their exact location. We introduce a mixed integer non-linear formulation of the problem and propose two linear reformulations. The latter involve sample average approximations and closest assignment constraints, and in one of the formulations we use scenario aggregation to reduce its size.
We solve the formulations with a general-purpose solver using a standard Benders decomposition method. 
Based on extensive computational results and a realistic case study, we find that the problem can be solved efficiently.
However, the level of uncertainty in the instances determines which approach is the most efficient. We use an entropy measure to capture the level of uncertainty\rev{, which} can be computed prior to solving. Furthermore, the results highlight the value of accurate demand modeling, as customer preferences have an important impact on the solutions and associated costs. 

\end{abstract}

\section{Introduction} \label{sec:intro}

The sustained increase in e-commerce transactions and online sales has promoted a rapid and sharp increase in demand for last-mile deliveries \citep{Rohmer2020CIRRELT202011AG}.  Moreover, customers are demanding fast deliveries with high service levels and expect flexible shipping options \citep[][]{doi:10.1287/mnsc.2020.01274}. This imposes increasing pressure on the distribution network of logistics providers (LP\rev{s}). Moreover, it may negatively impact the environment by increasing emissions due to failed and repeated deliveries, as well as contributing to increased congestion \citep{Molin_Kosicki_vanDuin_2022}. Alternative delivery options, in particular, collection and delivery points (CDP), arise as means of reducing costs and emissions \citep[][]{YUEN20181}. 

Using CDPs as an alternative to home delivery is not uncommon and its adoption varies around the world. For example, in Singapore, 5.5\% of online shoppers have adopted parcel lockers as a delivery option \citep{choo2016impact}, while this figure is 20\% in France and 10\% in Germany \citep{MORGANTI2014178}. The main impacts of CDPs are (i) reducing delivery costs thanks to improved consolidation and (ii) reduced return levels due to better first-delivery rates. Some customers may also perceive it as an improved service level. 

In this work, we address a strategic decision-making problem faced by a LP seeking to reduce costs by locating different types of CDPs to divert demand from its home delivery service. The problem crucially depends on the demand for services via CDPs, that is, it depends on customers' preferences. In this context, two sources of uncertainty arise. \emph{First}, at a strategic planning level, the LP does not have perfect knowledge about its future customers and their precise locations. \emph{Second}, the LP has imperfect information about customers' preferences. We assume that the LP has access to historical data capturing past customers' characteristics and preferences, and that such data can be used to estimate a random utility maximization (RUM) model to predict demand \citep{McFa81}. We propose to model our problem with a mixed integer non-linear formulation that integrates a RUM model of customer behavior. 

Our work is related to the growing body of literature on facility location problems where user's behavior is described with RUM models, also known as \emph{choice-based facility location} problems. Similar to several studies \citep[e.g.,][]{PACHECOPANEQUE202126, PinzonEtAl24, HAASE2014689, legault2023modelfree}, we use a sample average approximation (SAA) to derive a linear formulation for any additive RUM model. Closest to our work are the studies addressing facility location with fixed capacities  \citep{haase_management_2013, PinzonEtAl24}.
Also relevant to our work is the literature on location routing problems (LRPs), and, in particular, studies where not all customers are assigned to a delivery route and where they can get serviced by visiting a facility \citep{ARSLAN20211, STENGER2012702}. 
In this context, only a few studies consider customers' preferences \citep{JANJEVIC201937, GUERREROLORENTE2020106433} and, unlike our model, they consider the preferences to be perfectly known (i.e., described by a deterministic demand model). Similar to some works on LRPs \citep{JANJEVIC201937, GUERREROLORENTE2020106433}, 
we do not solve a routing problem due to its computational difficulties. Instead, we use routing cost estimations, in our case based on average routing costs from historical data.

We contribute to the literature \rev{on choice-based facility location problems} by addressing a capacitated CDP location problem where there is uncertainty about customers' locations and preferences. In contrast to the closely related work of \cite{haase_management_2013} and \cite{PinzonEtAl24}, we tackle two challenging aspects \rev{in this class of problems}: (i) potentially conflicting objectives between the LP and its customers, and (ii) strict assignments of customers to facilities. \rev{The importance of guaranteeing strict assignments is twofold: to ensure preference-maximizing allocations and to enable accurate routing cost estimations.} The main methodological contribution of our work is the introduction of formulations that can be \rev{efficiently} solved using a standard Benders decomposition method implemented in general-purpose solvers. This is achieved through a specific formulation of the follower (customer) problem in a bilevel formulation based on SAA. We demonstrate how this follower problem leads to the well-known \emph{closest assignment constraints} in our single-level linear reformulation. Furthermore, we present an extensive computational study analyzing the computing times and the effects of uncertainty \rev{in the demand information (preferences and locations)} on the solutions. Similar to \cite{legault2023modelfree}, we introduce an entropy measure and show how computational performance varies with the level of uncertainty in the instances. \rev{In addition, this work highlights the critical role of customer location accuracy for solution quality and problem tractability}. Finally, we showcase the practical applicability of our work by detailing a realistic case study of our industrial partner. 

The remainder of this paper is structured as follows: Section~\ref{sec:litreview} reviews related work and Section~\ref{sec:pbdescription} describes the problem we address in this work. Section~\ref{sec:model_form} introduces the non-linear mathematical model that we propose, and Section~\ref{sec:Iexample} details an illustrative example. In Section~\ref{sec:Sim_models}, we derive mixed integer linear programming (MILP) reformulations and discuss the closest assignment constraints. Section~\ref{sec:results} reports the computational study and Section~\ref{sec:study_case} presents the case study. Section~\ref{sec:conclusion} concludes.

\section{Literature Review} \label{sec:litreview}

This section reviews the literature related to our work. It mainly concerns facility location problems integrating choice models. However, we briefly comment on the literature related to LRPs as our location problem is impacted by routing costs, although we do not model that problem in detail. Finally, we overview the main contributions of our work.

There is a growing body of literature on facility location problems including RUM models, also known as choice-based facility location problems \citep{LIN2024}. We distinguish the works that address the problem using closed-form probabilities \citep[e.g.,][]{LJUBIC201846, MAI2020874, krohn_preventive_2021, duong_joint_2023} from those that address them by simulating utilities \citep[e.g.,][]{haase_management_2013, PACHECOPANEQUE202126, legault2023modelfree, doi:10.1287/ijoc.2022.0185, PinzonEtAl24}. Our work resides within the latter group but, unlike most studies, we consider capacitated facilities. \rev{Whereas capacitated facility location problems are well studied \citep[e.g.,][]{FiscEtAl17,FischEtAl16},} we are only aware of two works -- \cite{haase_management_2013} and \cite{PinzonEtAl24} -- that treat \rev{strict capacities in general (non-closed form) choice-based facility location problems} 
\rev{\citep[note that, e.g.,][also consider capacities, but these are enforced within individual scenarios rather than for all scenarios, which makes their approach unsuitable for our problem]{PACHECOPANEQUE202126}.}
Our work differs from \cite{haase_management_2013} in three main aspects. \emph{First}, in their case, the open facilities must satisfy all of the demand, while we consider an outside option. \emph{Second}, we consider a bilevel formulation because the objective of the LP and the customers are not necessarily aligned. In their case, the facility location objective function is aligned with the objective of the users as both maximize utility\rev{, which does not require a bilevel modelling approach}. \emph{Third}, they do not make strict customer assignments to facilities, whereas we do. 

\cite{PinzonEtAl24} consider a pricing and location problem where the random utilities are associated with price and service levels. \rev{However,  these utilities do not depend on facility location decisions}. 
Moreover, similar to \cite{haase_management_2013}, \cite{PinzonEtAl24} do not make strict customer assignments to facilities. \rev{In addition, unlike \cite{PinzonEtAl24}, our approaches rely on exact methods rather than heuristics to solve the formulation resulting from SAA.}

Other works have addressed competitive facility location problems \citep[see, e.g.,][]{beresnev_capacitated_2016, beresnev_exact_2018, CASASRAMIREZ2018369} and capture problems \citep[see e.g.,][]{berman_locating_2007} with capacities. These works use deterministic patronizing rules instead of RUM models.

In regard to applications, the problem of locating parcel lockers is receiving increasing attention 
\citep[see,][for a review]{Rohmer2020CIRRELT202011AG}. We identify two classes of related work. The \emph{first} class assumes perfect knowledge of customer preferences \citep[][]{doi:10.1080/00207543.2017.1395490, LUO2022105677}. That is, the demand models are deterministic. On the contrary, the \emph{second} class involves models that do not assume perfect knowledge about demand and does not use closed-form stochastic demand models \citep[e.g.,][]{lin_last-mile_2020, LIN2022102541}. Our work 
defines a new class where the simulation approach opens the formulation to any type of RUM model.

Our work seeks to locate facilities so as to reduce estimated routing costs. However, we do not explicitly model the routing problem as is done in the LRP literature. This is a rich area of research \citep[see e.g.,][]{https://doi.org/10.1111/itor.12950, https://doi.org/10.1111/itor.12032, PRODHON20141} covering several problem variants \citep[e.g.][]{DREXL2015283}. Thus, several studies on LRP consider the location of depots that constitute supply points with the objective of serving every customer with delivery tours starting at one of the depots \citep[see e.g.,] []{NAGY2007649, doi:10.1287/trsc.2017.0746, ARSLAN20211}. There are also problem variants closer to our setting, such as location or routing problems \citep[see][]{ARSLAN20211} and LRP with subcontracting options \citep[see][]{STENGER2012702}. 
However, there is limited work on problem settings with parcel lockers \citep[see e.g.,][]{su8080828, JANJEVIC201937, GUERREROLORENTE2020106433}. Only \cite{JANJEVIC201937} and \cite{GUERREROLORENTE2020106433} consider customer's preferences for the pick-up locations and both studies assume perfect knowledge of customers' preferences. Moreover, they focus on challenges associated with approximating routing costs. Instead, we focus on modeling demand uncertainty, and we estimate routing cost using historical data.
  
Table~\ref{tab:related_work} summarizes the key aspects of our work in comparison to the literature \rev{thereby motivating our contributions}. We consider a CDP location problem with capacitated facilities and a stochastic demand model where the customers' utility-maximizing objective may be conflicting with the facility location objective. Uncertainty arises because, at the strategic planning horizon, customer locations are not perfectly known. Moreover, we assume imperfect knowledge of customers' preferences which we capture with an additive RUM model without making any specific assumptions about its type. \rev{As highlighted in the table, our work differs from the literature on choice-based facility location problems (i.e., a choice model that is not deterministic) by solving a model that considers capacitated alternatives \emph{and} conflicting objectives \emph{while} ensuring strict assignments. In addition, we consider uncertainty in customers' locations.} The next section introduces the problem in greater detail.

\setlength{\arrayrulewidth}{0.25pt} 

\setlength{\tabcolsep}{2pt}
\begin{landscape}
\begin{table}[]
\footnotesize
\centering
\begin{tabular}{l|cc|ccc|ccc|}
\cline{2-9}
& \multicolumn{2}{|c}{Application} 
& \multicolumn{3}{|c|}{Choice Model} &
\multicolumn{3}{|c|}{Optimization Model}\\       
\hline
  \multicolumn{1}{|l|}{Reference} &
  \begin{tabular}[c]{@{}l@{}}Parcel \\ Lockers\end{tabular} &
  \begin{tabular}[c]{@{}l@{}} \rev{Uncertainty in} \\ \rev{Customers' Locations} \end{tabular} &
  \multicolumn{1}{c}{\begin{tabular}[c]{@{}l@{}}Type \\  of Model\end{tabular}} &
  \begin{tabular}[c]{@{}l@{}}\rev{Endogenous} \\ \rev{Variables} \end{tabular} &
  \begin{tabular}[c]{@{}c@{}} \rev{Capacitated}  \\ \rev{Alternatives}\end{tabular} &
  \begin{tabular}[c]{@{}c@{}}Strict \\ Assignments\end{tabular} &
  \begin{tabular}[c]{@{}l@{}}Bilevel \\ program\end{tabular} &
  \multicolumn{1}{c|}{\begin{tabular}[c]{@{}l@{}} Conflicting\\ \multicolumn{1}{c}{objectives}\end{tabular}} \\ 
  \hline
\multicolumn{1}{|l|}{\cite{berman_locating_2007}} & &
   & D 
   & \rev{-} &
  \checkmark &  
  \checkmark& &
   \\ \hdashline[0.5pt/1pt] 
\multicolumn{1}{|l|}{\cite{beresnev_capacitated_2016}} & &
   & D
   & \rev{-} &
  \checkmark & 
  \checkmark& 
  \checkmark & \checkmark \\ \hdashline[0.5pt/1pt] 
\multicolumn{1}{|l|}{\cite{beresnev_exact_2018}} & &
   & D
   & \rev{-} &
  \checkmark & 
  \checkmark&
  \checkmark  & \checkmark \\ \hdashline[0.5pt/1pt] 
\multicolumn{1}{|l|}{\cite{CASASRAMIREZ2018369}} & &
   & D
   & \rev{-} &
   & &
  \checkmark & \checkmark \\ \hdashline[0.5pt/1pt] 
\multicolumn{1}{|l|}{\cite{doi:10.1080/00207543.2017.1395490}} &
  \checkmark  & & D
   & \rev{Location} &
   & & & \\ \hdashline[0.5pt/1pt] 
\multicolumn{1}{|l|}{\cite{duong_joint_2023}} & &
   &
    MNL & \rev{Location} &
   & & &
   \\ \hdashline[0.5pt/1pt] 
\multicolumn{1}{|l|}{\cite{GUERREROLORENTE2020106433}} &
  \checkmark & & D
   & \rev{-} &
   & & &
   \\ \hdashline[0.5pt/1pt] 
\multicolumn{1}{|l|}{\cite{haase_management_2013}} & &
   &
  Any & \rev{Location} &
  \checkmark &  & &
   \\ \hdashline[0.5pt/1pt] 
\multicolumn{1}{|l|}{\cite{JANJEVIC201937}} & 
  \checkmark & & D
   & \rev{-} &   
   & & &
   \\ \hdashline[0.5pt/1pt] 
\multicolumn{1}{|l|}{\cite{krohn_preventive_2021}} & &
   &
  MNL & \rev{Location} &
   & & &
   \\ \hdashline[0.5pt/1pt] 
\multicolumn{1}{|l|}{\cite{doi:10.1287/ijoc.2022.0185}} & &
   &
  Any & \rev{Location} &
   & &
  \checkmark & \\ \hdashline[0.5pt/1pt] 
\multicolumn{1}{|l|}{\cite{lin_last-mile_2020}} & 
  \checkmark & &
  MNL & \rev{Location} &
   & & & \\ \hdashline[0.5pt/1pt] 
\multicolumn{1}{|l|}{\cite{LIN2022102541}} & 
  \checkmark & & TLM & \rev{Location}
   & & & & \\ \hdashline[0.5pt/1pt] 
\multicolumn{1}{|l|}{\cite{LJUBIC201846}} & &
   &
  MNL & \rev{Location} &
   & & &
   \\ \hdashline[0.5pt/1pt] 
\multicolumn{1}{|l|}{\cite{LUO2022105677}} & 
  \checkmark & & D
   & \rev{-} &
   & & &
   \\ \hdashline[0.5pt/1pt] 
\multicolumn{1}{|l|}{\cite{MAI2020874}} & &
   &
  MNL & \rev{Location} &
   & & &
   \\ \hdashline[0.5pt/1pt] 
\multicolumn{1}{|l|}{\cite{MENDEZVOGEL2023834}} & &
   &
  PBL & \rev{Location and Price} &
   & &
  \checkmark & \checkmark\\ \hdashline[0.5pt/1pt] 
\multicolumn{1}{|l|}{\cite{legault2023modelfree}} & &
   &
  Any & \rev{Location} &
   & &
  \checkmark & \\ \hdashline[0.5pt/1pt] 
  \multicolumn{1}{|l|}{\cite{PinzonEtAl24}} & &
   &
  Any & \rev{Price and Service} & 
   & &
  \checkmark & \checkmark \\ \hdashline[0.5pt/1pt] 
\multicolumn{1}{|l|}{\cite{su8080828}} & 
  \checkmark & &
   & &
   & & &
   \\ \hline \hline
\multicolumn{1}{|l|}{Our work} &
  \checkmark & \checkmark &
  Any & \rev{Location} &
  \checkmark & 
  \checkmark&
  \checkmark & \checkmark\\
\hline 
\multicolumn{7}{l}{D: Deterministic} \\
\multicolumn{7}{l}{MNL: Multinomial Logit Model} \\
\multicolumn{7}{l}{TLM: Threshold Luce Model} \\
\multicolumn{7}{l}{PBL: Partially Binary Logit Model} \\
\end{tabular}
\caption{\rev{Summary of related work}}
\label{tab:related_work}
\end{table}
\end{landscape}
\setlength{\tabcolsep}{6pt}

\section{Problem Description} \label{sec:pbdescription}

This section details our strategic decision-making problem that consists of determining the locations for CDPs -- drop-off and pick-up facilities -- to minimize routing, operating, and implementing costs resulting from such a decision over a one-year planning period. We begin by introducing the types of facilities involved, then discuss the specific characteristics of the demand, and conclude with the impact of routing costs.   

The LP provides a parcel delivery service directly to its customers' home locations (standard home delivery). As a complement, it seeks to implement additional delivery services using CDPs. These CDPs are alternative locations where customers can pick up or drop off their parcels \citep[see, e.g.][]{JANJEVIC201937, Rohmer2020CIRRELT202011AG}. We distinguish two main types of CDPs: pick-up (P) points, often integrated into stores or local shops, and automated parcel lockers, which do not require staff assistance during pick-up or drop-off. The LP considers two types of parcel lockers: Regular (A) and modular (M). Regular parcel lockers have fixed capacity, whereas in modular parcel lockers, the capacity can increase by multiples of a base capacity. Let $D$ denote the set of candidate locations for CDPs. Also, let $D^\text{P}$, $D^\text{A}$ and $D^\text{M}$ denote the corresponding sets for pick-up points P, and parcel lockers A and M, respectively. Each $d \in D$ has coordinates $\theta_d=(\theta_d^{\text{x}}, \theta_d^{\text{y}}) \in \mathbb{R} ^2$. Moreover, $d \in D^\text{P} \cup D^\text{A}$ has a fixed capacity $u_d$, whereas deciding the capacity of $d \in D^{\text{M}}$ is part of our problem. 

The LP does not have precise knowledge about future customers and deliveries. Customers are, therefore, aggregated based on their geographical location and preferences. More precisely, customers are classified into categories $k \in K$ according to delivery profiles. These delivery profiles are based on customers' characteristics and preferences captured in historical data. Furthermore, we consider two levels of geographical aggregation. CDPs attract demand from customers in a zone $z$ within a set of zones $Z$. To model demand, \rev{given the lack of precise information on customer locations}, each zone is further partitioned into subzones $a \in A_z$ \rev{capturing uncertainty in the spatial distribution of demand}. Customers are distributed across subzones, each having an expected yearly number of parcels $b_{zak},~z\in Z, a\in A_z, k\in K$ with a total number of parcels within a zone $b_z=\sum_{a\in A_z} \sum_{k\in K} b_{zak}$.

The demand that each CDP attracts is influenced by the customers' preferences for this option relative to the standard home delivery service. 
The distance from the customer's location to a CDP is an important explanatory factor. Given the geographical aggregation of customers, we consider the distance from the centroid of the corresponding subzone $a$, denoted by $\theta_a \in \mathbb{R}^2$. Other factors can impact the customers’ preferences and these are not all known to the LP. This motivates the use of a stochastic demand model in our case. Note that, as the demand captured by a CDP cannot exceed its capacity, excess demand is reallocated to the home delivery service.

Implementing CDPs implies variable (i.e., volume-dependent) and fixed costs $f_d > 0~,d \in D^\text{A}\cup D^\text{M}$. Note that fixed costs increase with capacity for $d \in D^\text{M}$. Moreover, typically $f_d=0,~d \in D^\text{P}$, but these pick-up points require a minimum demand level $\underline{b}_d$ from the LP to respect business agreements with the stores.

Variable costs are mainly routing costs to serve the customers.  The LP estimates the current unit routing costs per parcel based on historical delivery and routing operations for each zone, $c_{z}^\text{R},~z\in Z$. They encompass both transit costs, incurred when the vehicle moves between a terminal to the corresponding zones and within the zone, as well as delivering-operating costs related to loading, sorting, and moving parcels from the vehicle to the corresponding delivery addresses.

The use of CDPs impacts the estimated routing costs as fewer customers require home delivery services, and the consolidated demand from captured customers requires less transit and delivering-operating costs. In this problem, explicit routing decisions cannot be made as the LP uses customer data aggregated over subzones. However, the LP considers that there is an impact on the delivering-operating costs due to parcel consolidation when delivering to a CDP instead of to different home addresses. Hence, the LP quantifies the affected routing (AR) costs of customers in a zone $z$ served by a CDP $d$ by assuming a reduction in the delivering-operating costs, regardless of the transit costs.  We denote these costs by $c_{dz}^\text{AR}$. Note that in this case $c_{dz}^\text{AR} \leq c_{z}^\text{R}$.

In summary, the LP seeks to optimally locate CDPs to minimize total expected costs encompassing estimated routing costs and CDP's fixed and operational costs. Since these costs depend on the demand for CDPs, the LP must consider customers' preferences for CDPs, their capacity as well as the preferences for home delivery services.

\section{A Mixed Integer Non-Linear Formulation}\label{sec:model_form}

This section introduces a non-linear formulation for our problem. It is the foundation of the two linear reformulations that are derived and discussed in Section~\ref{sec:Sim_models}. We begin by presenting the main modeling concepts and then we focus on the mathematical formulation itself. 

We denote by $x_d$ the binary decision of implementing or not a CDP at $d \in D$ (in vector notation $\boldsymbol{x} \in \{0,1\}^{|D|}$). We introduce demand models at the two levels of geographical aggregation (zone and subzone) that depend on $\boldsymbol{x}$. 

The share of customers in $z$ choosing delivery service $d$ is denoted $\psi_{dz}(\boldsymbol{x})$. Since the captured demand is constrained by the CDP's capacity, we define the effective share of demand captured by a CDP $d$ in zone $z$ and denoted it with the auxiliary variables $p_{dz}$. Note that the latter implicitly depends on $\boldsymbol{x}$ since $p_{dz} \leq \psi_{dz}(\boldsymbol{x})$.

Customer information is available at the subzone level. We model demand for CDPs and home delivery services with a RUM model and we denote the probability of customers in category $k$ and subzone $a$ choosing $d$ by $\psi_{dak}(\boldsymbol{x})$. Accordingly, the demand shares, at zone level, are defined by 
\begin{equation}
    \psi_{dz}(\boldsymbol{x})=\sum_{a \in A_{z}}\sum_{k \in K} \frac{b_{zak}}{b_z}\rev{\rho_{dak}(\boldsymbol{x})}, \label{prob_dzkq}
\end{equation}
with
\begin{equation}\label{eq:ch_prob_a}
    \rev{\rho_{dak}(\boldsymbol{x})} = \text{Prob} \left (d = \underset{d' \in D_x \cup  \{0\}}{\arg \max} \{U_{d'ak}(\boldsymbol{\varepsilon_{d'}})\} \right ).
\end{equation}
In (\ref{eq:ch_prob_a}), $U_{d'ak}(\boldsymbol{\varepsilon_{d'}})$ is an additive random utility function capturing the attractiveness of  
a CDP $d'$ for customers of category $k$ in subzone $a$. It is defined over the set of open CDPs $D_{\boldsymbol{x}}$ and home delivery service (indexed with $0$). For simplicity, we henceforth write $U_{dak}$ instead of $U_{d'ak}$ and define it as 
\begin{equation} 
    U_{dak} = \alpha^{\text{1}}_{k} f^\text{1}(\theta_d, \theta_a)+\alpha^{\text{2}}_{k}f^\text{2}(E_d) +\varepsilon_{dak}, \label{assp:general_utility} 
\end{equation} 
where $f^\text{1}(\theta_d, \theta_a)$ denotes a function of the distance from the customers located in subzone $a \in A_z$ to the location of the CDP $d$, and $f^\text{2}(E_d)$ is a function of other features, exogenous to our model, that may explain customer preferences. 
Parameters $\alpha^{\text{1}}_{k}$ attached to CDP alternatives are strictly negative and reflect that the probability of choosing a CDP decreases as distance increases. For home delivery service, $\alpha^{\text{1}}_{k}=0$. Parameters $\alpha^{\text{2}}_{k} \geq 0$ reflect the relative importance of other attributes for customer category $k$. Terms $\varepsilon_{dak}$ are continuous random variables with support $\Xi$. Specific distributional assumptions lead to different types of RUM models.  For example, assuming that $\varepsilon_{dak}$ are independently and identically distributed \rev{(i.i.d.)} extreme value type I leads to the well-known logit model. 

\rev{A couple of notes regarding~\eqref{assp:general_utility} are in order. First, in practice, parameter values should be estimated from data. The estimation of RUM models is well-studied \citep[e.g.,][]{Train09} and typically performed independently of the downstream optimization problem. Whereas there is a growing body of literature on integrating estimation and optimization \citep{SADANA2025271}, no existing methodology addresses our type of problem, where decision-dependent uncertainty appears in the constraints of the mathematical program. In this work, we do not focus on model estimation. Instead, we fix parameter values to assess the impact of the RUM model on solutions through sensitivity analyses. Second, as we further explain in Section~\ref{sec:Sim_models}, we restrict the utility function to be linear in the decision variables to obtain a linear single-level formulation of the problem. No specific assumptions are needed for the part of the utility function that depends only on exogenous features.}

Modular parcel lockers $d \in D^\text{M}$ have capacities that increase by multiples of a base capacity. Accordingly, we define a set of capacity levels $L_d~=\{u^\text{1}_d, u^\text{2}_d, \ldots, u^\text{L}_d \}$, each with a corresponding fixed cost denoted by $f^l_d$. Moreover, we introduce binary decision variables $r^l_d$ that equal one if level $l$ is selected for CDP $d$, and zero otherwise.

The mixed integer non-linear formulation of our problem is: 
\newtheorem{assumption}{Assumption}
\begin{subequations}\label{model1}
    \begin{equation}\label{enlm:fo} 
        \begin{split}
            \text{min } & \sum_{z \in Z}  c_{z}^\text{R}b_zp_{0z}+  \sum_{z \in Z}\sum_{d \in D} c_{dz}^\text{AR}b_zp_{dz} + \sum_{d \in D^\text{M}} \sum_{l \in L_d} f^l_d r^l_d +\sum_{d \in D^\text{A}} f_d x_d 
        \end{split}
    \end{equation}
    \begin{align}
        \text{s.t.~~} &\sum_{l \in L_d} r^l_{d} = x_d, & &  d \in D^\text{M} \label{enlm:1} \\
         &\sum_{z \in Z}  b_{z} p_{dz} \leq \sum_{l \in L_d} u^l_d r^l_d, & & d \in D^\text{M} \label{enlm:2} \\
         &\sum_{z \in Z}  b_{z} p_{dz} \leq u_d x_d, & & d \in D^\text{P} \cup D^\text{A} \label{enlm:2.5} \\
         &\sum_{z \in Z} b_z p_{dz} \geq \underline{b}_d x_d, & & d \in D^P \label{enlm:2.8} \\
         & p_{dz} \leq \psi_{dz}(\boldsymbol{x}), & & z \in Z, d \in D \label{enlm:3} \\
         & p_{0z} \geq \psi_{0z}(\boldsymbol{x}), & & z \in Z \label{enlm:p_0} \\
         & p_{0z}+\sum_{d \in D} p_{dz} = 1, & & z \in Z \label{enlm:4} \\
         & x_d \in \{0,1\}, & & d \in D \label{enlm:dom_x} \\
         &p_{dz} \geq 0, & & z \in Z, d \in D \cup 0 \label{enlm:dom_p} \\
         & r^l_d \in \{0,1\}, & &  d \in D^
        \text{M}, l \in L_d. \label{enlm:dom_r}
    \end{align}
\end{subequations}
Objective function~\eqref{enlm:fo} minimizes the total estimated routing costs, including the costs of implementing and operating CDPs. It has four components capturing the (i) total estimated routing costs for home-delivery service, (ii)  total estimated routing and operating costs related to demand allocated to open CDPs, as well as total fixed costs of implementing (iii) modular and (iv) regular parcel lockers, respectively.
Constraints~\eqref{enlm:1} impose that only one of the available capacity levels can be chosen for modular parcel lockers. Constraints~\eqref{enlm:2}~and~\eqref{enlm:2.5} ensure that the expected demand captured by a CDP does not exceed its capacity, whereas Constraints~\eqref{enlm:2.8} impose a minimum demand level for pickup points $d \in D^\text{P}$. Constraints~\eqref{enlm:3} impose that the effective fraction of demand captured by a CDP cannot exceed the demand share. In contrast,  Constraints~\eqref{enlm:p_0} state that the effective fraction of demand captured by the home-delivery service must be at least equal to the demand share. We include these constraints for the sake of clarity but note that they are, in fact, redundant.
 
Constraints~\eqref{enlm:4} ensure that all the demand is captured by either a CDPs or the home delivery service. 

We note that model~(\ref{model1}) is non-linear due to Constraints~\eqref{enlm:3}, in view of the definitions of $\psi_{dz}(\boldsymbol{x})$ in~\eqref{prob_dzkq} and of the individual  probabilities in~\eqref{eq:ch_prob_a}. Indeed, even the simplest RUM model -- multinomial logit -- is non-linear. For more advanced RUM models, such as the logit mixtures \citep{Train09}, \eqref{eq:ch_prob_a} does not even have a closed form whence simulation is required to evaluate probabilities. In this situation, one way to solve model~(\ref{model1}), is to introduce simulated utilities that are maximized as part of a follower problem in a bilevel formulation \citep[e.g.][]{PACHECOPANEQUE202126, doi:10.1287/ijoc.2022.0185,PinzonEtAl24}. In our case, the bilevel structure is important, as the RUM objective of the user may not be aligned with the cost minimization objective of the LP. We describe such reformulations in Section~\ref{sec:Sim_models} after the following example designed to illustrate the impact of capacities and costs on the location decisions and captured demand. Finally, we refer to Tables~\ref{tab:art2:notation_part1} and~\ref{tab:art2:notation_part2} in Appendix~\ref{art2:Main_notation} for a comprehensive listing of the notation.

\section{Illustrative Example}\label{sec:Iexample}

This section illustrates the impact of locating CDPs in a region divided into two zones, where customers are served by the LP's standard home delivery service. We discuss how capacity constraints affect the demand captured by these facilities and highlight the balance the LP must achieve between capturing demand and minimizing costs.

\begin{figure}[htbp]
    \includegraphics[width=5cm]{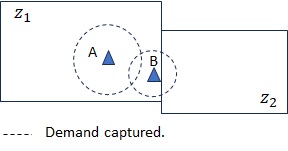}
     \centering
    \caption{Zones and CDPs}
    \label{fig:Iexample}
\end{figure}

We consider two zones $z_1$ and $z_2$ and two candidate locations for CDPs $\{A,B\}$, as depicted in Figure~\ref{fig:Iexample}.
\begin{table}[htbp]
\small
\centering
\begin{tabular}{lrrr}
\hline
$\psi(\boldsymbol{x})$ & \multicolumn{1}{l}{$x_A=1$} & \multicolumn{1}{l}{$x_B=1$} & \multicolumn{1}{l}{$x_A=x_B=1$} \\
\hline
$\psi_{Az_1}$ & 0.35 & 0    & 0.30  \\
$\psi_{Bz_1}$ & 0    & 0.20  & 0.15 \\
$\psi_{0z_1}$ & 0.65    & 0.80 & 0.55 \\
\hline
$\psi_{Az_2}$ & 0    & 0    & 0    \\
$\psi_{Bz_2}$ & 0    & 0.05 & 0.05 \\
$\psi_{0z_2}$ & 1    & 0.95 & 0.95 \\
\hline
\end{tabular}
\caption{Choice probabilities for each CDP location decision}
\label{tab:Iprob}
\end{table}

The demand in each zone is given $b_{z_1}=b_{z_2}=500$, respectively. The unit costs of home delivery service in each zone are $c_{z_1}^\text{R}=5$ and $c_{z_2}^\text{R}=6$. For the sake of simplicity, we assume a single customer category and no partition into subzones. The capacity for the CDP in each location is given by $u_A=90, u_B=130$, and the fixed costs are $f_A=170, f_B=200$. 

Figure~\ref{fig:Iexample} shows the relative location of CDPs and highlights that, from a distance perspective, they are more attractive for customers in zone $z_1$. We fix the choice probabilities as reported in Table~\ref{tab:Iprob} for different location decisions $x_d$.  Recall that the effective fraction captured by a CDP can be lower than the corresponding probabilities $\psi$, due to the capacity constraints. 

The estimated costs of serving customers captured by the CDPs are $c_{Az_1}=c_{Az_2}=3.25$ and $c_{Bz_1}=c_{Bz_2}=3.5$. If none of the CDPs are open, the total estimated routing cost is \$5,500 ($5\cdot 500 + 6 \cdot 500$). This corresponds to serving all the customers with home-delivery services, thus $p_{0z_1}=p_{0z_2}=1$. However, the optimal solution is $x_A=0, x_B=1$, with a total cost of \$5,488. In this solution, CDP B captures 20\% of the demand from zone A and 5\% from zone B. Even though CDP A is more attractive and less expensive than CDP $B$, its capacity limits the amount of demand captured, impacting the potential cost reduction. We note that the optimal solution is sensitive to the capacities of CDPs. For example, if the capacity of CDP A increases by at least three units, then the optimal solution changes to $x_A=1, x_B=0.$  Table~\ref{tab:Iresults} summarizes the results.

\begin{table}[htbp]
\centering
\footnotesize
\begin{tabular}{lllcc}
\hline
                      & $x_A=1$                & $x_B=1$                 & \multicolumn{1}{c}{$x_A=x_B=1$} & \multicolumn{1}{c}{$x_A=x_B=0$} \\
                      \hline
Total cost           & 5,513                 & 5,488                  & 5,538                            & 5,500                            \\[2pt]
Captured demand zone $z_1$ & \multicolumn{1}{c}{$90$} & \multicolumn{1}{c}{$100$} & $165$                             & 0                               \\
Effective fraction zone $z_1$ & \multicolumn{1}{c}{18\%} & \multicolumn{1}{c}{$20
\%$} & $33\%$                             & 0                               \\[2pt]
Captured demand zone $z_2$ & \multicolumn{1}{c}{0}  & \multicolumn{1}{c}{$25$}  & $25$ & 0  \\
Effective fraction zone $z_2$ & \multicolumn{1}{c}{0}  & \multicolumn{1}{c}{5\%}  & $5\%$    
& 0                              
\end{tabular}
\caption{Results for location decisions $x_d$}
\label{tab:Iresults}
\end{table}

\normalsize

The results in Table~\ref{tab:Iresults} also highlight the impact of the conflicting objectives between customers and the LP. 
We observe that a solution capturing more demand (more customers allocated to their preferred choice) is not necessarily the best for the LP. Indeed, when both CDPs are open, they capture more demand. However, this solution is more costly than using only home delivery services for all the customers.

Recall the rule that demand is allocated to home delivery service in case of excess demand. When both CDPs are open, the effective fraction of demand captured follows the LP's allocation rule. This fraction results in a lower bound on the effective fraction that both CDPs can capture when customers can be reallocated to the CDP next in preference. To show this, let us consider that the excess of demand captured by CDP $A$ is reallocated to CDP $B$ until its maximum capacity is reached. Then, the total captured demand from zone $z_1$ would be $195$. That corresponds to $39\%$ ($>33\%$) of the demand in zone $z_1$ and represents a lower cost (\$5,493) than the one displayed in Table~\ref{tab:Iresults} for the same solution ($x_A=x_B=1$). However, in this case, the corresponding total cost with this allocation rule is still higher than the cost obtained by the optimal solution under the allocation assumption in this work.

\section{Simulation-based Formulations}\label{sec:Sim_models}

This section presents our approaches to deal with the non-linear program~\eqref{model1} introduced in Section~\ref{sec:model_form}. They are bilevel MILP formulations based on SAA of customer choice probabilities. We begin this section by discussing the foundation of the SAA formulation and two modeling alternatives for the follower’s problem. Then, we focus on the modeling alternative leading to the single-level formulation with the smallest number of variables and constraints. We introduce in Section~\ref{sec:CAC} the single-level reformulation, and then, in Section~\ref{sec:SAAA}, an equivalent reformulation based on scenario aggregation. 

We adopt an existing approach that consists in simulating realizations of random utilities that are linear in decision variables  \citep[e.g.][]{haase_management_2013, PACHECOPANEQUE202126}. \rev{Several works discuss the convergence of SAA in this context; see, for example, Section~3.2 in \cite{PACHECOPANEQUE202126} and Section~3.3 in \cite{legault2023modelfree}. These rely on the asymptotic convergence properties of SAA \citep{KleywegtEtAl02,Kim2015} and the ability to draw i.i.d. samples from the utility distributions (variance reduction techniques may also be applied).}  We draw an \rev{i.i.d.} sample set $S$ of scenarios $\{\boldsymbol{\varepsilon_s}\}_{s \in S}$ from the distribution of the random terms $\varepsilon_{dak}$ in~\eqref{assp:general_utility}. We estimate $\rev{\rho_{dak}(\boldsymbol{x})}$ with an SAA through
\begin{subequations}
\small
\begin{align}
    \rev{\hat{\rho}_{dak}(\boldsymbol{x})} & =\frac{1}{|S|}\sum_{s \in S} w_{daks}, & d\in D, a \in A, k \in K \label{SAA_prob}  \\
    w_{daks}& = \begin{cases}\label{SAA_newvar}
        1, & d =\underset{d' \in D_x \cup  \{0\}}{\arg \max} \{U_{d'aks}\}  \\
        0, &\text{otherwise.}  
    \end{cases} & d \in D, a \in A, k \in K, s \in S
\end{align}
\end{subequations}
In~\eqref{SAA_newvar}, $w_{daks}$ equals one if delivery service $d$ (a CDP or home-delivery service) has the maximum utility in scenario $s$ for customers in category $k$ and subzone $a$. \rev{Note that only one $d \in D_x \cup 0$,  for each $a, k, s$,  attains value 1, since $U_{dkas}$ are drawn from continuous distributions.} These $w_{daks}$ are decision variables in the follower problem of our bilevel formulation. We consider two formulations. For given $a$, $k$, and $s$, and (leader) decision variables $\boldsymbol{x}$, a first, straightforward formulation is
\begin{subequations}\label{af:P1}
             \begin{align}
                {\arg \max} & \sum_{d \in D} U_{daks}x_d w_{daks} +U_{0aks}w_{0aks} \label{eq:P1.of} \\
                \text{s.t. } & \sum_{d \in D \cup \{0\}} w_{daks} =1, \label{eq:uniqued} \\  
                & w_{daks} \in \{0,1\}, & d \in D \cup \{0\},
            \end{align}
        \end{subequations}
and a second, alternative formulation is 
\begin{subequations}\label{af:P2}
        \begin{align}
        {\arg \max} &\sum_{d \in D \cup \{0\}} U_{daks}w_{daks} \label{eq:P2.of} \\
        \text{s.t. } & \sum_{d \in D \cup \{0\}} w_{daks} =1, \label{eq:uniqueedP2} \\
        & w_{daks} \leq x_d, & d \in D \label{eq:P2.w} \\ 
                & w_{daks} \in \{0,1\}, & d \in D \cup \{0\}.
        \end{align}
    \end{subequations}
    
We note that the coefficient matrices in these two formulations are totally unimodular. Therefore, we can relax the integrality constraints on variables $w_{daks}$ and still obtain binary solutions. We also note that without loss of generality, we assume that $U_{daks}>0$. Indeed, in the case of negative utilities, we can add a constant (e.g., $\hat{U}_{daks}=U_{daks}-\delta$, where $\delta<\min_{d' \in D \cup \{0\}} U_{d'aks}$) without affecting the solutions. Based on this assumption, and the fact that $U_{daks}$ are sampled from continuous distributions (the probability that any two values are equal is zero), the optimal solutions to the formulations are unique. In addition, in the following proposition, we show that the optimal solutions of the two formulations are equal.

\newtheorem{Prop}{Proposition}
    \begin{Prop}\label{Prop:sameobj}
        For given $a$, $k$, $s$, let $\hat{\boldsymbol{w}}_{aks}$ and $\tilde{\boldsymbol{w}}_{aks}$ be the vectors representing the corresponding optimal solutions for \eqref{af:P1} and \eqref{af:P2} given a set of open facilities $D^{*} \subset D$ defined by $\boldsymbol{x}$. Then $\hat{\boldsymbol{w}}_{aks}=\tilde{\boldsymbol{w}}_{aks}$.
    \end{Prop}
    \begin{proof}
        By inspection of~\eqref{af:P1} and~\eqref{af:P2}, we see that any level of utility reachable in the objective of~\eqref{af:P1} can also be reached in~\eqref{af:P2} by adjusting $\tilde{\boldsymbol{w}}_{aks}$ and conversely that any level of utility reachable in~\eqref{af:P2} can also be reached in~\eqref{af:P1} by adjusting $\hat{\boldsymbol{w}}_{aks}$. Hence the sets of utility levels reachable in~\eqref{af:P1} and~\eqref{af:P2} are identical.
       For simplicity, we henceforth omit indices $a,k$, and $s$. We have $x_d=1, d \in D^{*}$. Let us assume an optimal solution to~\eqref{af:P1}: $\hat{w}_{d^{'}}=1, d^{'} \in D^{*} \cup \{0\}$, and $\hat{w}_d=0, d \neq d^{'} \in D^{*} \cup \{0\}$.
       This implies that $U_{d^{'}} >U_{d}, d \neq d^{'} \in D^* \cup \{0\}$. 
       Let us assume a different optimal solution to~\eqref{af:P2}: $\tilde{w}_{d^{''}}=1, {d^{''} \neq d^{'} } \in D^{*} \cup \{0\}$, and $\tilde{w}_d=0, d \neq {d^{''}} \in D^{*} \cup \{0\}$.
       This implies that $U_{d^{''}} > U_d, d \neq {d^{''}} \in D^{*} \cup \{0\}$. However, this contradicts the assumption that 
       $\hat{\boldsymbol{w}}$ is optimal for \eqref{af:P1}.
       Hence, both optimal solutions must be equal.    
    \end{proof}

In the following section, we introduce a single-level formulation using~\eqref{af:P2} where a subset of the constraints are closest assignment constraints (CAC). In Section~\ref{sec:SAAA}, we introduce another single-level formulation that is also based on~\eqref{af:P2}, but that uses an aggregation of scenarios inspired by \cite{legault2023modelfree}.

\subsection{CAC Formulation}\label{sec:CAC}
Next, we introduce a linear single-level formulation of \eqref{model1}. It is derived by replacing $\psi_{dz}(\boldsymbol{x})$ \rev{and $\psi_{0z}(\boldsymbol{x})$} in Constraints~\eqref{enlm:3} \rev{and~\eqref{enlm:4} } by their SAA using \rev{the alternative formulation for \eqref{SAA_newvar} given by }~\eqref{af:P2}. This results in a bilevel formulation. Then, using the strong duality theorem, we derive the following single-level and linear reformulation (henceforth referred to as formulation C): 

\footnotesize
\begin{subequations}\label{model:Sim_CAC}
\begin{equation}
    \begin{split}
        \text{min } & \sum_{z \in Z}  c_{z}^\text{R}b_zp_{0z}+  \sum_{z \in Z}\sum_{d \in D} c_{dz}^\text{AR}b_zp_{dz} + \sum_{d \in D^\text{M}} \sum_{l \in L_d} f^l_d r^l_d +\sum_{d \in D^\text{A}} f_d x_d  \label{saam:fo}
    \end{split}
\end{equation}
\begin{align}
    \text{s.t: } & \eqref{enlm:1}- \eqref{enlm:2.8}, \eqref{enlm:4} - \eqref{enlm:dom_r}, \nonumber \\
     & \sum_{g \in D \cup \{0\}} e_{gdaks} w_{gaks}  \geq x_d, & & d \in D, a \in A, k\in K,  s \in S \label{saam:max_ut_1} \\
      &\sum_{g \in D \cup \{0\}} e_{g0aks} w_{gaks} \geq 1, & &  a \in A, k \in K, s \in S \label{saam:max_ut_2} \\
     &\sum_{d \in D \cup\{0\}} w_{daks}  = 1, & & a \in A, k \in K, s \in S \label{saam:del_serv_assign} \\
     &w_{daks} \leq x_d, & & d \in D, a \in A, k \in K, s  \in S \label{saam:bound_w} \\
     & p_{dz} \leq \frac{1}{|S|}\sum_{s \in S}\sum_{a \in A_{z}}\sum_{k \in K}\frac{b_{zak}}{b_z} w_{daks}, & & d \in D, z \in Z \label{saam:bound_prob} \\
     & p_{0z} \geq \frac{1}{|S|}\sum_{s \in S}\sum_{a \in A_{z}}\sum_{k \in K}\frac{b_{zak}}{b_z} w_{0aks}, & &  z \in Z \label{saam:bound_prob_0} \\
     &0 \leq w_{daks} \leq 1, & &  d \in D\cup \{0\},  a \in A, k \in K, s  \in S. \label{saam:dom_w}
\end{align}
\normalsize

\end{subequations}
\normalsize
The objective function~\eqref{saam:fo}, sets of constraints as well as restrictions on the decision variables are the same as in~\eqref{model1}, 
except that  Constraints~\eqref{saam:max_ut_1}-~\eqref{saam:dom_w} and decision variables $w_{dajs}$ are substituted for Constraints~\eqref{enlm:3}.

As we show below, Constraints~\eqref{saam:max_ut_1} - \eqref{saam:bound_w} result from applying the strong duality theorem to \eqref{af:P2}. Specifically, Constraints~\eqref{saam:max_ut_1} and~\eqref{saam:max_ut_2} enforce the assignment to the CDP having the highest utility. In such constraints, the parameter $e_{gdaks}=1$ if $U_{gaks} \geq U_{daks}, d \in D$, and $0$ otherwise. Similarly, $e_{g0aks}=1$ if $U_{gaks} \geq U_{0aks}$, and 0 otherwise. Moreover, Constraints~\eqref{saam:del_serv_assign} ensure that customers in subzones are assigned to a single CDP or to the home delivery service. Right-hand-sides of Constraints~\eqref{saam:bound_prob} and Constraints~\eqref{saam:bound_prob_0} act as the SAA of $\psi_{dz}(\boldsymbol{x})$ and $\psi_{0z}(\boldsymbol{x})$, respecively. 

Note that Constraints~\eqref{saam:max_ut_1} and~\eqref{saam:max_ut_2} correspond to one of the most widely used CAC \citep[see Constraints~($\mathcal{CC}$) on page 50 in][]{RePEc:eee:ejores:v:219:y:2012:i:1:p:49-58}. However, in our case, we seek utility-maximizing assignments as opposed to the classic distance-minimizing assignments. 

In the following, we show that the CAC-constraints given in \eqref{saam:max_ut_1}, \eqref{saam:max_ut_2} with Constraints~\eqref{saam:del_serv_assign} and \eqref{saam:bound_w}, are derived from the follower problem \eqref{af:P2}. 
\rev{To this end, we introduce the dual formulations of~\eqref{af:P1} and~\eqref{af:P2}, then we rewrite them using linear programming strong duality, and prove that this reformulation is equivalent to the CAC-constraints~\eqref{saam:max_ut_1} - \eqref{saam:bound_w}}

\rev{Although this section focuses on the single-level reformulation resulting from the follower’s problem \eqref{af:P2}, we use the dual of ~\eqref{af:P1} to replace the optimal solution of the dual of \eqref{af:P2}. Both formulations yield the same optimal solution as noted in Proposition~\ref{Prop:sameobj}, however, the dual of \eqref{af:P1} admits an analytical expression, which facilitates the derivation of the CAC-constraints.}

\rev{The duals for~\eqref{af:P1} and~\eqref{af:P2} are given by \eqref{af:D.1} and \eqref{af:D.2}}
\begin{subequations}\label{af:D.1}
    \begin{align}
   \min \ & \alpha_{aks} \label{af:D1:of} \\
    \text{s.t. } & \alpha_{aks} \geq  U_{daks} x_d, & d \in D \label{eq:d1.alp1} \\
    & \alpha_{aks} \geq U_{0aks}. \label{eq:D1.alp2}
\end{align}
\end{subequations}
and
\begin{subequations}\label{af:D.2}
\begin{align}
    \min \ & \gamma_{aks} +\sum_{d \in D} \beta_{daks}x_d \label{af:D2:of} \\
    \text{s.t. } & \gamma_{aks} +\beta_{daks} \geq  U_{daks}, & d \in D \label{eq:D2.alp1} \\
    & \gamma_{aks} \geq U_{0aks}, \label{eq:D2.alp2} \\
    & \beta_{daks} \geq 0, & d \in D.  \label{eq:D2.beta}
\end{align}
\end{subequations}
        
Using strong duality, we rewrite~\eqref{af:P1} and \eqref{af:P2} as the following constraint satisfaction problems:
\begin{align}
    \sum_{d \in D} U_{daks} x_d w_{daks} +U_{0aks}w_{0aks} \geq \alpha_{aks} \label{Sl-P1} \\
    \eqref{eq:uniqued}, \eqref{saam:dom_w},\eqref{eq:d1.alp1}, \eqref{eq:D1.alp2} \nonumber
\end{align}
and
\begin{align}
    \sum_{d \in D \cup \{0\}} U_{daks} w_{daks} \geq \gamma_{aks} +\sum_{d \in D} \beta_{daks} x_d\label{Sl-P2} \\
    \eqref{eq:uniqueedP2}, \eqref{eq:P2.w}, \eqref{saam:dom_w}, \eqref{eq:D2.alp1}, \eqref{eq:D2.alp2}, \eqref{eq:D2.beta}, \nonumber
\end{align}
respectively. We use these when proving the result in the following proposition.
\begin{Prop}\label{CAC-derive}
    For each $a,k$ and $s$, \rev{the constraint satisfaction problem~\eqref{Sl-P2} is equivalent} to CAC-constraints~\eqref{saam:max_ut_1}, \eqref{saam:max_ut_2}, \rev{together with Constraints}  \eqref{saam:del_serv_assign}, \eqref{saam:bound_w}.

\end{Prop}
\begin{proof}
    First, note that Constraints~\eqref{saam:del_serv_assign} and \eqref{saam:bound_w} are the same as Constraints~\eqref{eq:uniqueedP2} and \eqref{eq:P2.w}. Thus, we focus on showing that Constraints~\eqref {Sl-P2}, \eqref{eq:D2.alp1}, \eqref{eq:D2.alp2} and \eqref{eq:D2.beta} yield the CAC Constraints~\eqref{saam:max_ut_1} and \eqref{saam:max_ut_2}.
 
    According to the strong duality theorem, Constraints~\eqref {Sl-P2} are only satisfied (at the equality) at the optimum. Hence:
    \begin{align}
        \sum_{g \in D \cup \{0\}} U_{gaks} w_{gaks} & \geq \gamma^{*}_{aks} +\sum_{d \in D} \beta^{*}_{daks} \label{Proof.P3.sd2} \\
         & \geq \alpha_{aks}^{*} \label{Proof.P3.P2} \\
        & \geq \max\{U_{daks} x_d, d \in D; U_{0aks}\} \label{Proof.P3.SolD1}
    \end{align}

    In inequality ~\eqref{Proof.P3.P2} we use Proposition~\ref{Prop:sameobj} to replace the optimal objective value of \eqref{af:D.2} with the optimal objective value of \eqref{af:D.1}, as from there it is evident that $\alpha_{aks}^{*}= \max\{U_{daks} x_d, d \in D; U_{0aks} \}$.         
    We replace \eqref{Proof.P3.sd2}-\eqref{Proof.P3.SolD1} with the following two sets of constraints:
    \begin{align}
         \sum_{g \in D \cup \{0\}} U_{gaks} w_{gaks} & \geq U_{daks} x_d, d \in D \label{proof:cac_D} \\
         \sum_{g \in D \cup \{0\}} U_{gaks} w_{gaks} & \geq U_{0aks}. \label{proof:cac_0}
    \end{align}
   Note that only for $g \in D \cup \{0\}$ where $U_{gaks} \geq U_{daks}$, can $w_{gaks}$ satisfy Constraints~\eqref{proof:cac_D}. Also, only those with $g \in D \cup \{0\}$ where $U_{gaks} \geq U_{0aks}$, $w_{gaks}$ can satisfy Constraint~\eqref{proof:cac_0}.  
   By using $e_{gdaks}=1$ for $U_{gaks} \geq U_{daks}$; 0 otherwise, and $e_{g0aks}=1$ for $U_{gaks} \geq U_{0aks}$; 0 otherwise, Constraints~\eqref{proof:cac_D} and \eqref{proof:cac_0} can be written as:
                   
     \begin{align}
         \sum_{g \in D \cup \{0\}}  e_{gdaks}  w_{gaks} & \geq  x_d, d \in D \label{proof:cacr_d} \\
         \sum_{g \in D \cup \{0\}} e_{g0aks} w_{gaks} & \geq 1, \label{proof:cacr_0}
    \end{align}
    which corresponds to Constraints~\eqref{saam:max_ut_1} and \eqref{saam:max_ut_2}, \rev{the CAC-constraints}.
\end{proof}

We close this subsection with a few remarks on the solution method. Problem~\eqref{model:Sim_CAC} can be solved using a standard Benders decomposition method implemented in a general-purpose solver. We identify two strategies to partition the variables: The \emph{first} one is \emph{automated partition} (AP) where integer variables are assigned to the master problem and continuous variables to the subproblems. A \emph{second} strategy, that we call \emph{user-partition} (UP), includes also the continuous variables $p_{dz}$ in the master problem. By doing so, we can decompose the subproblems by zones. However, the resulting subproblems are feasibility subproblems because none of the variables appear in the objective function. Hence, the UP only generates feasibility cuts. We use C-AP and C-UP, to respectively refer to the approach using AP or UP to solve model~\eqref{model:Sim_CAC}. 

\subsection{CAC with Scenario Aggregation }\label{sec:SAAA}

In this section, we describe a formulation that is based on the CAC formulation of Section~\ref{sec:CAC} and the SAA with \emph{aggregation} model (SAAA) introduced in \cite{legault2023modelfree}. In this work, they achieve strong computational performance by aggregating scenarios that lead to the same solution. Moreover, they propose a partial Benders decomposition where they only explicitly account for scenarios that contribute sufficiently to the objective function. Their solution approach, however, cannot be directly applied to our model because our objective function is not submodular, and we have capacity constraints. In the following, we explain how we can nevertheless define scenario aggregation in a way that is similar to theirs. 

Let $(a,k,s)$ for $a \in A_z$, $k \in K$ and $s \in S$ denote a triplet. For each $(a,k,s)$ we have a vector of utilities $\boldsymbol{U_{aks}}~=~[U_{1aks}, U_{2aks} , \ldots , U_{Daks}, U_{0aks} ]$. For each vector $\boldsymbol{U_{aks}}$, we associate a vector $\Pi_{aks}$ denoting the ranking of the corresponding utilities. For example, for $|D|=2$, let us assume the following subset of triplets $(a,k,s)$: $\{(1,1,1),(1,1,2), (1,2,1),(3,3,2)\}$, with the corresponding utilities: $\{[4,7,3], [2,5,1], [3,6,1], [2,5,0]\}$. Accordingly, $\Pi_{111}=\Pi_{112}=\Pi_{121}=\Pi_{332}=[2,3,1]$. Consider for instance the utility vector $[4,7,3]$, in this case the utilities of 4, 7 and 3, have rankings 2, 3 and 1, respectively.
All of the utility vectors above have the same ranking. 

For the purposes of aggregation, we define a pattern  $q$ that identifies a subset of triplets $(a,k,s)$ such that the associated vector of ranked utilities, $\Pi_{aks}$, are equal. Let $Q$ denote the set of unique patterns. Moreover, let $\nu_{qak}$ denote the fraction of scenarios $s \in S$ in which customers in subzone $a \in A_z$ in category $k \in K$ appear in pattern $q \in Q$. \rev{Algorithm~\ref{alg:patterns} provides the steps to derive the patterns and the associated parameters.} 

\begin{algorithm}
    \caption{\rev{Compute  patterns $q \in Q$}}\label{alg:patterns}
    \begin{algorithmic}[1]
        \Require $U_{daks}, \text{for each } d,a,k,s$
        \ForAll {$a$, $k$, $s$}
	   \State $\Pi \gets \text{argsort} (U_{daks}, d \in D \cup 0)$ \Comment{Get the ranking of locations}
	   \State $q \gets \text{Encode}(\Pi)$ \Comment{Encode ranking into pattern}
	   \State $\overline{U}[q] \gets \Pi$ \Comment{Store ranking}
	   \State $\nu_{ak}[q] \gets \text{get}(\nu_{ak} [q],0) + 1$ \Comment{Update pattern frequency count}
       \EndFor
        \State \Return $\nu, \overline{U}$
    \end{algorithmic}
\end{algorithm}
\rev{In Algorithm~\ref{alg:patterns}, in Line~2, the ranking order of locations $d \in D \cup 0$ is based on the value of the simulated utilities. Note that this ranking is computed for each $a, k, s$. In Line 3, $
\text{Encode}(\cdot)$ may be constructed using several encoding strategies (e.g., concatenation of elements or hashing). This pattern is unique for a given $\Pi$, but it may be shared across different $a, k, s$ if the same ranking $\Pi$ occurs multiple times. In Line 5, we update a dictionary $\nu_{ak}$ which keeps track of the frequency of each pattern. If a pattern $q$  appears for the first time, the function $\text{get}(\cdot)$ returns the default value $0$; otherwise, it retrieves the current number of cases where pattern $q$ has appeared.}

Compared to the previous formulation, we slightly change the definition of our decision variables $w_{dq}$ that now equal 1 if the observed pattern $q$ is assigned to CDP $d$, and 0 otherwise. In addition, we use $e_{gdq}=1$, if \rev{$\overline{U}_{gq} \geq \overline{U}_{dq}$} and 0, otherwise, for $g,d \in D \cup \{0\}$. This leads to a formulation with aggregation (referred to as AC) where the difference with respect to~\eqref{model:Sim_CAC} resides in the use of patterns $q \in Q$ instead of scenarios $s \in S$:

\footnotesize
\begin{subequations}\label{model:Sim_CAC-A}
\begin{equation}
    \begin{split}
        \text{min } & \sum_{z \in Z}  c_{z}^\text{R} b_z p_{0z}+  \sum_{z \in Z}\sum_{d \in D} c_{dz}^\text{AR}b_zp_{dz} + \sum_{d \in D^\text{M}} \sum_{l \in L_d} f^l_d r^l_d +\sum_{d \in D^\text{A}} f_d x_d  \label{pat:fo}
    \end{split}
\end{equation}
\begin{align}
    \text{s.t: } & \eqref{enlm:1}- \eqref{enlm:2.8}, \eqref{enlm:4} - \eqref{enlm:dom_r}, \nonumber \\
     & \sum_{g \in D \cup \{0\}} e_{gdq} w_{gq}  \geq x_d, & & \forall d \in D, q \in Q, \label{pat:max_ut_1} \\
      &\sum_{g \in D \cup \{0\}} e_{g0q} w_{gq} \geq 1, & & \forall q \in Q \label{pat:max_ut_2} \\
     &\sum_{d \in D \cup \{0\}} w_{dq} = 1, & &\forall q \in Q \label{pat:del_serv_assign} \\
     &w_{dq} \leq x_d, & &\forall d \in D, q \in Q \label{pat:bound_w} \\
     & p_{dz} \leq \sum_{q \in Q}\sum_{a \in A_{z}}\sum_{k \in K}\frac{b_{zak}}{b_z} \nu_{qak} w_{dq}, & & \forall d \in D, z \in Z \label{pat:bound_prob} \\
     & p_{0z} \geq \sum_{\rev{q \in Q}}\sum_{a \in A_{z}}\sum_{k \in K}\frac{b_{zak}}{b_z} \nu_{qak} w_{0q}, & &  z \in Z \label{pat:bound_prob_0} \\
     &0 \leq w_{dq} \leq 1, & & \forall d \in D\cup \{0\}, q \in Q. \label{pat:dom_w}
\end{align}

\end{subequations}

\normalsize

As with model~\eqref{model:Sim_CAC}, we solve~\eqref{model:Sim_CAC-A} using a standard Benders decomposition method with AP partition (\rev{henceforth referred to as} AC-AP). We note that UP is not suitable due to the linking Constraints~\eqref{pat:bound_prob} and~\eqref{pat:bound_prob_0}. Indeed, the patterns $q \in Q$ may aggregate customers from different zones, hindering the decomposition over zones.

\section{Computational Results}\label{sec:results}
This section describes an extensive computational study designed to achieve two main objectives: \emph{First}, to analyze how the levels of uncertainty in the demand model and customer locations impact the solutions. \emph{Second}, to assess the computational performance with increasing problem size. In the following section, we describe the experimental setup, and then discuss the results related to each of the objectives in Sections~\ref{sec:Interpret} and~\ref{sec:Largsec}, respectively.

\subsection{Experimental Setup}\label{sec:SetInstances}
The experimental setup concerns the definition of the problem instances, model parameters, and performance metrics, as well as the characteristics of our infrastructure and the computing budget. We define the sets of instances by describing how we generate demand, and candidate CDP locations. For all experiments we consider a region represented by a $[0, 30] \times [0,30]$ plane, two customer categories ($|K|=2$), and a distribution center located at $[0,0]$ in the plane.

\paragraph{Demand Generation.} Recall that there are two sources of uncertainty, both pertaining to demand: The LP does not have perfect knowledge about the precise location of future customers or about their preferences. To assess the impact of the first source, we generate demand with a precise location, and then we assess the impact of aggregating customers into subzones of varying sizes. As for the second source, we divide customers into two categories. Namely, 25\% of the customers (drawn at random) are likelier to use CDPs than the other 75\% of the customers. The strength of their preferences is governed by parameters in the utility functions as described below. 

We pseudo-randomly sample $|J|$ customers' locations ($ \theta^\text{x}_j, \theta^\text{y}_j$) within the plane ($0 \leq \theta^\text{x}_j \leq 30$ and $0 \leq \theta^\text{y}_j \leq 30$). Note that the superscript indicates that these correspond to x and y-coordinates in the plane. We consider two different sampling protocols. In the first (U-R, shorthand for uniform random), we assume that customers are uniformly distributed over the region. In the second (UN-R, shorthand for uniform-normal random), we assume a higher concentration of customers at the center ($[15,15]$) of the region. In this case, 50\% of customers are uniformly distributed in the plane (same as U-R), and the other 50\% are distributed according to a normal distribution $\theta^\text{x}_j,\theta^\text{y}_j \sim N(15,3)$.

\paragraph{Demand Aggregation into Zones.} We partition the region into equally sized zones with the same number of subzones. Since we simulate the precise location of customers, we can aggregate all of those who are in a given zone. Hence, fewer zones implies more aggregation. 

In Section~\ref{sec:Interpret}, we use the same number of customers $|J|=4,000$, and zones $|Z|=4$, but vary the level of aggregation through the number of subzones $|A|=\{4, 8, 16, 64, 256\}$. In Section~\ref{sec:Largsec} where we focus on computational performance, we use $|J|=100,000$ customers and vary either the number of zones -- Z-class instances -- or the number of subzones -- A-class instances. In Z-class, we use $|Z|=\{32, 48, 64, 96\}$, and in A-class, we use $|A|=\{192, 256, 320, 400, 480, 576, 784, 1024\}$. 

Note that, given the pseudo-randomly generated precise customer locations, it is possible to have instances with subzones that have no customers, or only customers of a single category. Hence, the total combination of subzones and customer categories $|\tilde{A}| \leq|A||K|$. 

Next, we turn our attention to the generation of candidate locations for CDPs.

\paragraph{CDP Generation.} We pseudo-randomly sample $|D|$ CDP locations within the plane according to a uniform distribution. For each type of CDP, the capacity is the same. We consider the 3 types of CDPs, and $|D|/3$ locations for each type, where modular parcel lockers ($D^\text{M}$) have three capacity levels ($|L|=3$) and $L=\{u_d, 2u_d, 3 u_d\}$. The minimum demand level is $\underline{b}_d=0.5 u_d$ and fixed cost $f_d=2 u_d$. In Section~\ref{sec:Interpret}, aimed at analyzing the impact of uncertainty, we fix $|D|=15$ and $u_d=20$, whereas we generate more challenging instances in Section~\ref{sec:Largsec} -- referred to as D-class instances -- using $|D|=\{21,27,33,39,45\}$. We use $u_d =\{550, 825, 1650\}$ for Z, A and S-classes. For the D-class, we proportionally adjust the capacities such that the total capacity is the same for all instances.

\newcommand{\norm}[1]{\left\lVert#1\right\rVert}

In the following we describe model parameters, starting with how we compute variable costs, followed by the RUM model.

\paragraph{Variable Costs.} We estimate routing costs in zone $z\in Z$ with
$c_{z}^\text{R}=0.1 \times \norm{\theta_z}$, where $\norm{\cdot}$ denotes the $L_2$-norm. Note that this corresponds to the Euclidean distance between the distribution center located at [0,0] and the centroid of the zone.
The resulting routing cost are $c_{dz}^\text{AR} =0.90 \times c_z^\text{R}$, for $d \in D^\text{P}$, and $c_{dz}^\text{AR} =0.80 \times c_z^\text{R}$, for $d \in D \setminus D^\text{P}$.

\paragraph{RUM Model.} We postulate a logit model in all experiments, i.e., $\varepsilon_{daks}$ are \rev{i.i.d.} according to the Extreme Value type I distribution with scale parameter $\beta$. 
We use a random utility defined as $U_{daks} =$~$\alpha_k^1\big(|\theta^{\text{x}}_{a} -\theta^{\text{x}}_d|+|\theta^{\text{y}}_{a} -\theta^{\text{y}}_d|\big)+\varepsilon_{akds}$\rev{, for each scenario $s \in S$, where the size of $S$ varies according to the class of instances as we describe below}. Here, again, the superscripts denote the x- and y-coordinates in the plane. We use $\alpha_2^1 = 5 \alpha_1^1$ so that customers in category $k=2$ are five times more sensitive to the distance than customers in category $k=1$. \rev{This captures the idea that some customers are less willing to travel to CDPs to collect their parcels}. Hence, the values of $\beta$ and $\alpha_1^1$ determine the level of uncertainty and customers' preferences. 

In Section~\ref{sec:Interpret}, we use three configurations of these parameters  $(\alpha_1^1, \beta)=\{(-0.01, 1), (-0.1,1), (-0.1, 0.25)\}$, whereas in Section~\ref{sec:Largsec} we use a subset of two configurations $(\alpha_1^1, \beta)=\{(-0.1, 1), (-0.1, 0.25)\}$.

\paragraph{Instances and Scenarios.} Our formulations are based on scenarios obtained by sampling values of the random terms. For each of the parameter configurations that we have described, we draw ten sets of scenarios. 
An instance refers to a problem defined by a set of scenarios, hence we have ten instances for each parameter setting.
Moreover, we consider different numbers of scenarios. In Section~\ref{sec:Interpret}, we use $|S|=50$, and we also keep it fixed to this number when solving the Z-, A- and D-classes of instances. To assess the impact on computing time, we consider a fourth class of instances, referred to as the S-class where we vary the number of scenarios $|S|=\{50, 100, 150, 200, 250, 300\}$.

\paragraph{Performance Metrics and Interpretability.} We report different metrics depending on the focus of the analysis, including computing time in seconds, the number of CDPs that are open in the optimal solution, the value of the objective function, the optimality gap for instances that are not solved to optimality within our time limit. For interpretability of the level of uncertainty associated with the demand model, we use the entropy measure introduced in \cite{legault2023modelfree}. The entropy is defined as $H= \sum_{q \in Q} - v_q \log (v_q)$, where $v_q$ is the fraction of scenarios where pattern $q$ is observed. We report entropy, or average entropy $\hat H$. Note that the maximum entropy (i.e. the entropy of a uniform distribution) may vary across instances. Therefore, we also report the maximum entropy and categorize instances as having low, medium or high entropy relative to the maximum.

\paragraph{Infrastructure and Computing Budget.} The experiments are performed on the Linux version of IBM ILOG CPLEX 22.1, running on an Intel  I7-9700K 8 cores CPU at 3.6GHz/4.9GHz,  with a time budget of one hour.

\subsection{Impact of Uncertainty}\label{sec:Interpret}

\begin{figure}[htbp]
    \centering
    \includegraphics[width=0.7\textwidth]{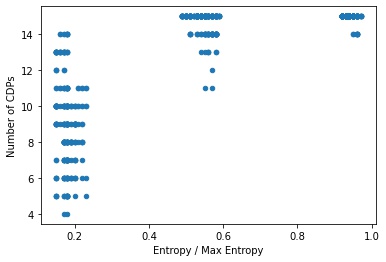}
    \caption{Number of open CDPs in the optimal solution for different values of the ratio between entropy and max entropy for each instance}
    \label{fig:opt_CDPS_entropy}
\end{figure}

\begin{figure}[htbp]
    \centering
    \subfloat[Medium entropy\label{fig:cdps_med_entropy}]{{\includegraphics[width=0.55\textwidth]{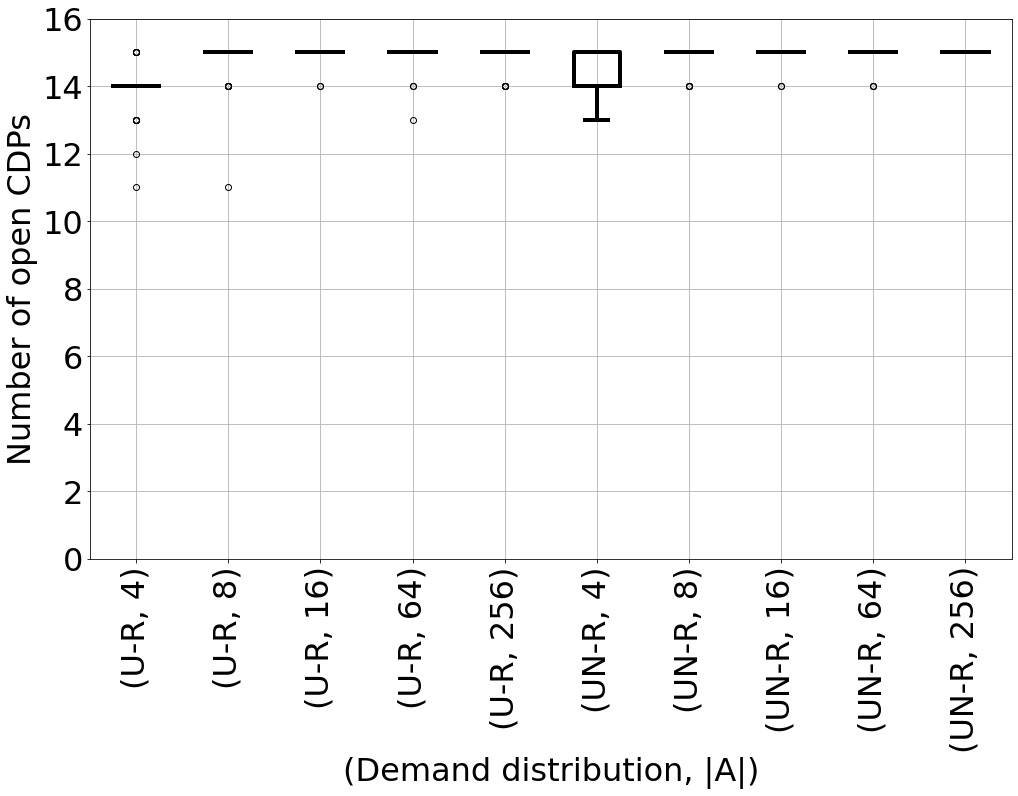}}}
    \subfloat[Low entropy \label{fig:cdps_low_entropy} ]{{\includegraphics[width=0.55\linewidth]{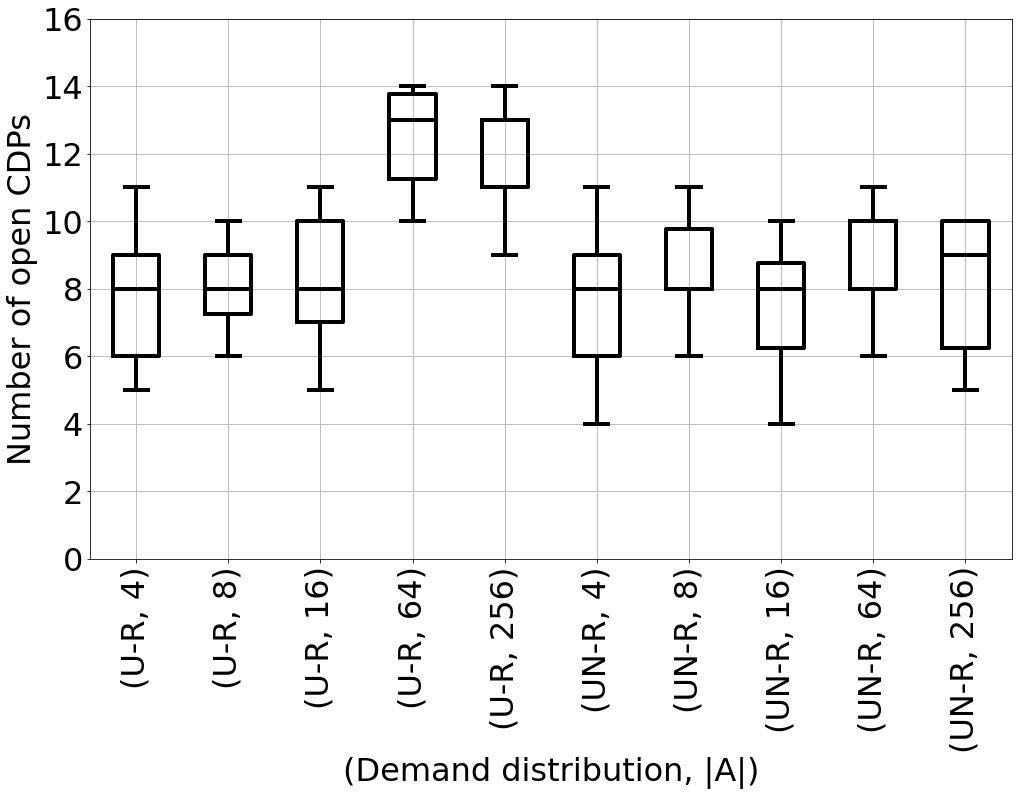}}}
    \caption{Distribution of the optimal number of CDPs}
    \label{fig:Cdps_dist}
\end{figure}

\begin{figure}[htbp]
    \centering
    \subfloat[Medium level entropy\label{fig:cost_med_entropy}]{{\includegraphics[width=0.55\textwidth]{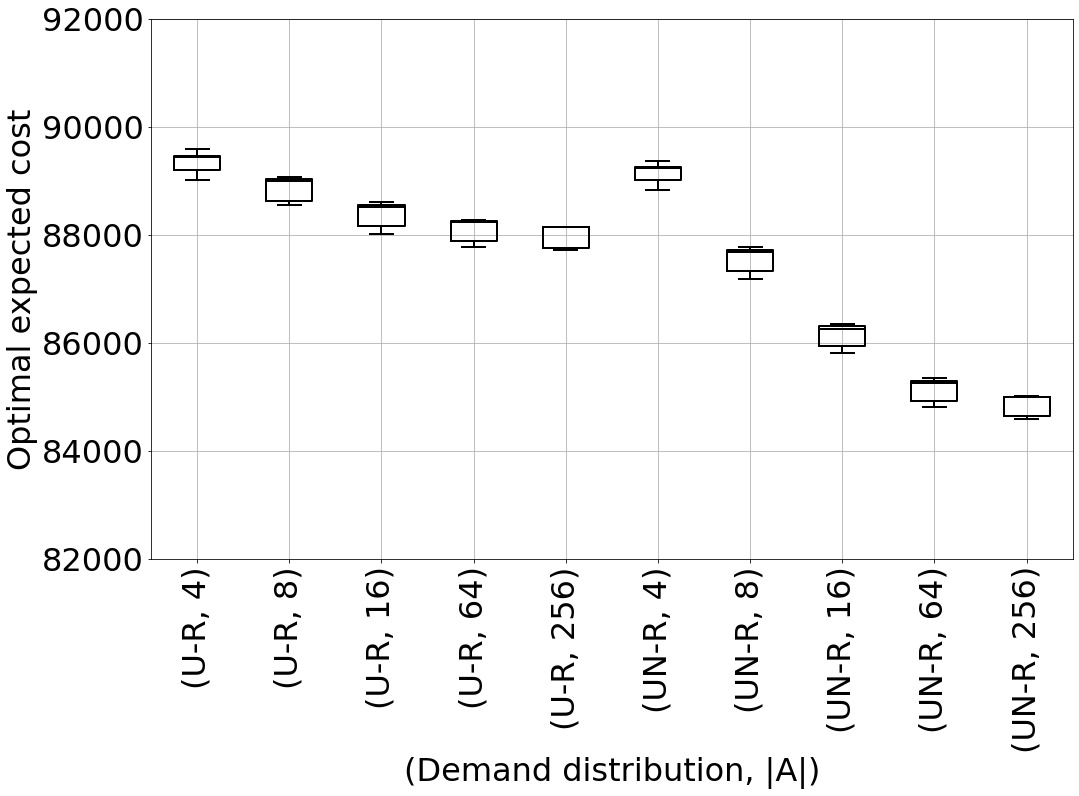}}}
    \subfloat[Low level entropy \label{fig:cost_low_entropy} ]{{\includegraphics[width=0.55\linewidth]{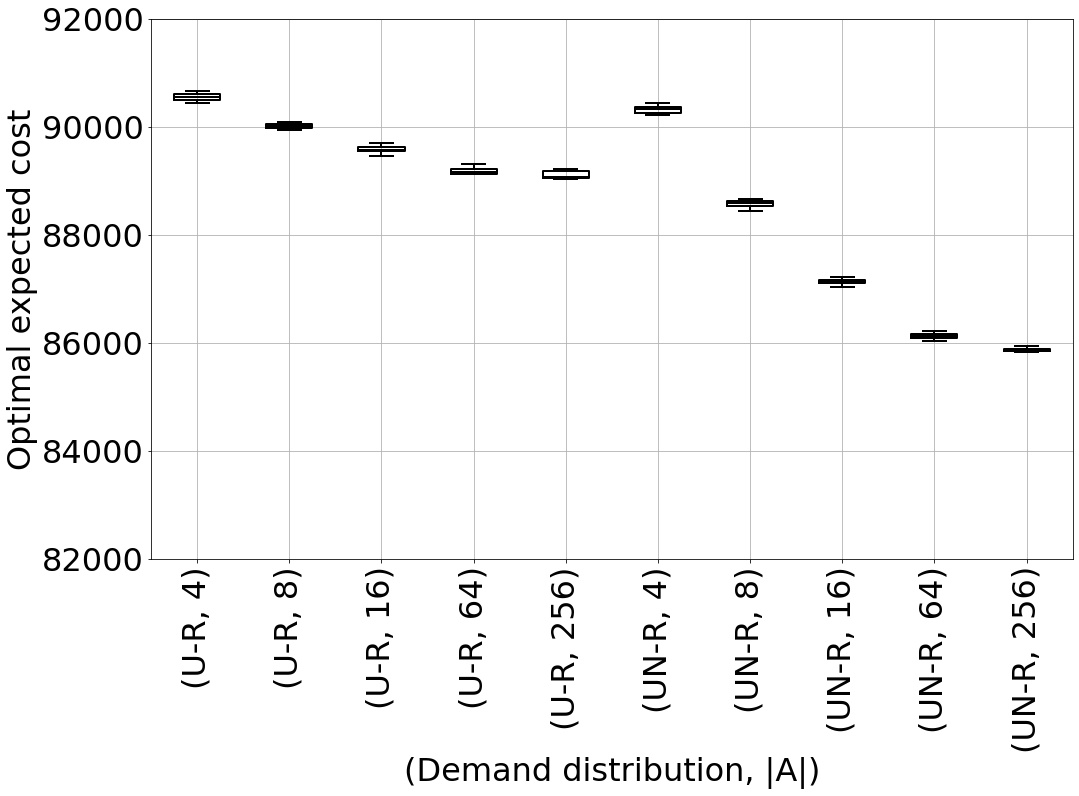}}}
    \caption{Distribution of the optimal expected cost}
    \label{fig:optimal_cost_Dist_A}
\end{figure}

In this section, we analyze two aspects that impact the level of uncertainty. Namely, (i) parameters of the RUM model, and (ii) geographical distribution of CDPs and customer aggregation into subzones. Recall that we keep the numbers of customers, zones and CDP fixed for the results in this section, while we use different configurations of the utility functions and different numbers of subzones. Figures~\ref{fig:opt_CDPS_entropy} --~\ref{fig:optimal_cost_Dist_A} display different views on the same set of instances.

Figure~\ref{fig:opt_CDPS_entropy} is a scatter plot showing, for each instance, the number of open CDPs in the optimal solution and the associated entropy level expressed as a ratio of the maximum entropy. All 15 CDPs are open in most high-entropy instances. 
High entropy implies that all CDPs have similar estimated probabilities to capture demand. In this case, the total reduction in estimated routing costs surpasses the fixed cost. Moreover, all CDPs have the same fixed cost.  
Therefore, the model chooses to open all CDPs. 

The figure clearly shows three clusters of instances: low, medium, and high entropy. 
Note that we include high-entropy configurations in these results for the sake of illustration. In practice, such configurations are of less interest as they imply demand models that do not provide any information compared to a uniform distribution. In the following we provide a more in-depth analysis of the low and medium entropy instances.

    Figures~\ref{fig:cdps_med_entropy} and~\ref{fig:cdps_low_entropy} display the distributions of the number of open CDPs in the optimal solutions as a function of parameter configurations. The latter is the geographical distribution of demand (uniform, U-R, and more concentrated to the center UN-R), and the number of subzones $|A|$. We note that all 15 CDPs are open in most of the medium entropy instances except those having the fewest number of subzones. We see a higher variance in the low-entropy instances and fewer open CDPs. The median (indicated with a line in each box) is higher for larger number of subzones, especially when demand is uniformly distributed in the region. 
    
    Now, turning to the optimal expected costs displayed similarly in Figures~\ref{fig:cost_med_entropy} and~\ref{fig:cost_low_entropy}, we note that they decrease as the number of subzones increases. Moreover, they are smaller when demand is more highly concentrated in the center of the region (UN-R compared to U-R).  \rev{As expected, the distribution of the optimal expected cost exhibits lower variance in low-entropy instances, as the deterministic components of the utility functions dominate the error terms. Conversely, the variance increases with higher entropy levels.}
    \rev{We report additional results in Appendix~\ref{apx:Add_results} where we vary the size of the set of scenarios, $|S|$. We note that, as expected, the variance decreases for larger numbers (we display results up to 400 scenarios) and the median value stabilizes. However, 50 scenarios result in relatively accurate estimates and using more scenarios does not change the analysis of the results that follows next.}

Recall that the model trades off reduction in variable routing costs with fixed costs associated with opening CDPs. The reduction of variable routing costs depends on the demand that a CDP can attract which, in turn, is governed by the utility functions. In the low-entropy instances, attraction strongly depends on distance, whereas this dependence is less strong in medium-entropy instances. 

We use the centroid of a subzone to compute distances. The distribution of demand in the region, along with the size of the subzones, therefore determine the quality of the distance approximation. Consider uniformly distributed demand in the region. Having only four subzones means that distance will be overestimated for a significant share of the demand, compared to having smaller (more) subzones. This issue is less pronounced when demand is concentrated to the center of the region (UN-R distribution). 

Given the above, less demand is attracted by CDPs, especially in low-entropy instances when subzones are relatively large (see Figure~\ref{fig:Cdps_frac_cap} in Appendix~\ref{apx:Add_results}). This, hence, explains why those instances have fewer open CDPs and higher costs. In other words, it is valuable to reduce uncertainty by trying to model demand as accurately as possible (utility functions and geographical distribution), and to solve the model with a larger number of subzones. Of course, from a practical point of view, the size of the subzones should reflect the level of uncertainty in the geographical distribution of demand. It does not make sense to have small subzones \rev{when} there is a high degree of uncertainty \rev{about} where the demand is located. As expected, areas with relatively high concentration of demand \rev{result in} lower costs. 


\subsection{Model and Algorithm Comparison}\label{sec:Largsec}

In this section we analyze the computing time when the size of the instances increases in various dimensions. We compare the three approaches C-AP, C-UP and AC-AP. As detailed in Section~\ref{sec:SetInstances}, we use four classes of instances for this purpose: S-, Z-, A- and D-classes.

Tables~\ref{tab:s-class_instance} --~\ref{tab:d-class_instance} report the average computing time (in seconds) for the four classes of instances. Recall that the average in each row is computed over ten instances. All instances are solved to optimality with only 10 exceptions (three configurations of the D-class, see footnotes in Table~\ref{tab:d-class_instance}). In addition to computing time, we report average entropy ($\hat{H}$), and average maximum entropy (avg max $H$). Each table is separated into two parts, the upper reports results for medium-entropy instances, and the lower for low-entropy instances. The shortest computing time in each row is highlighted in bold font.

A few findings clearly emerge from these results. The AC-AP is the best approach for low-entropy instances across all classes. This is expected, and consistent with the findings in \cite{legault2023modelfree}, as aggregation has the largest impact on model size for low-entropy instances. We note that AC-AP dominates the other approaches for each individual instance and not only on average (see Figure~\ref{fig:ratio_time_ac_ap} in Appendix~\ref{apx:Add_results}). In terms of absolute computing time, low-entropy instances are also the easiest to solve. This opens up the possibility to use large number of subzones and scenarios for improved modeling accuracy. 

Next, we focus on findings related to medium-entropy instances. First, we note that increasing the size of $Z$ marginally impacts the computing time as this dimension only affects continuous variables $p_{dz}$. The hardest instances arise when increasing $|D|$ as it governs the number of binary variables $x_d$. Instances where the capacity is not restrictive (larger $u_d$) are especially challenging. When capacities are restrictive, more demand is allocated to home delivery service, which makes the problem easier to solve. Computing time also increases with $|S|$ and $|A|$ but less strongly than increasing $|D|$.

Unlike for low-entropy instances, a single approach does not dominate the others. For the S-class, C-UP performs the best as $|S|$ increases. There is no clear winner for the A-class instances but the computing times are quite similar for the three approaches, except for the two largest configurations where AC-AP is significantly faster on average. On the hardest class of instances (D-class), C-AP performs the best overall. C-UP has a lower average computing time on the largest instances but solves only one out of ten instances to optimality within the time limit. For the two configurations where AC-AP has the lowest average computing time, C-AP has a similar performance. 

The entropy of an instance can be computed before it is solved. The results indicate that the solution approach should be selected based on this information. For low-entropy instances, the best approach is to use aggregation AC-AP. For medium (to high) entropy instances, C-AP consistently performs well, if not the best. Especially for the hardest instances with relatively high number of candidate CDP locations. However, if the objective is to improve modeling accuracy by increasing the number of scenarios and the number of subzones, then C-UP can achieve a strong performance. We close this section by noting that the approaches we propose based on reformulations using~\eqref{model:Sim_CAC} consistently display shorter computing times (sometimes by several orders of magnitude) compared to those that result from using~\eqref{af:P1}. We provide more details in Appendices~\ref{apx:DSLR} and~\ref{apx:Add_results}.

\begin{table}[htbp]
\centering
\footnotesize
\begin{tabular}{llll|rrr|}
\cline{5-7}
 &
   & &
  \multicolumn{1}{l|}{} &
  \multicolumn{3}{c|}{\textbf{Avg. computing time (sec)}} \\ \hline
\multicolumn{1}{|l}{\textbf{$\hat{H}$}} &
{avg max $H$} &
\textbf{$|S|$} &
\textbf{$u_d$} &
\textbf{C-AP} &
\textbf{C-UP} &
\textbf{AC-AP} \\ \hline
\multicolumn{1}{|c}{4.5} & 8.8 &50  & 550  & {\bf 2.0}   & 2.3  & 2.2   \\
\multicolumn{1}{|c}{4.5} & 8.8 &50  & 825  & {\bf 2.4}   & 2.6  & 2.6   \\
\multicolumn{1}{|c}{4.5} & 8.8 &50  & 1650 & {\bf 3.0}   & 4.4  & 3.3   \\
\multicolumn{1}{|c}{4.8} & 9.4 &100 & 550  & {\bf 5.4}   & 5.8  & 6.7   \\
\multicolumn{1}{|c}{4.8} & 9.4 &100 & 825  & {\bf 6.4}   & 6.7  & 7.5   \\
\multicolumn{1}{|c}{4.8} & 9.4 &100 & 1650 & {\bf 8.7}   & 9.5  & 10.1  \\
\multicolumn{1}{|c}{5.0} & 9.9 &150 & 550  & {\bf 10.6}  & {\bf 10.6} & 13.8  \\
\multicolumn{1}{|c}{5.0} & 9.9 &150 & 825  & 12.7  & {\bf 11.9} & 16.1  \\
\multicolumn{1}{|c}{5.0} & 9.9 &150 & 1650 & {\bf 17.1}  & 17.6 & 19.9  \\
\multicolumn{1}{|c}{5.1} & 10.1 &200 & 550  & 16.7  & {\bf 16.6} & 23.4  \\
\multicolumn{1}{|c}{5.1} & 10.1 &200 & 825  & 19.8  & {\bf 18.7} & 26.0  \\
\multicolumn{1}{|c}{5.1} & 10.1 &200 & 1650 & 27.8  & {\bf 27.6} & 33.3  \\
\multicolumn{1}{|c}{5.2} & 10.4 &250 & 550  & 25.3  & {\bf 24.1} & 35.7  \\
\multicolumn{1}{|c}{5.2} & 10.4 &250 & 825  & 30.0  & {\bf 26.8} & 39.9  \\
\multicolumn{1}{|c}{5.2} & 10.4 &250 & 1650 & 41.0  & {\bf 40.2} & 50.2  \\
\multicolumn{1}{|c}{5.3} & 10.5 &300 & 550  & 35.1  & {\bf 33.9} & 49.5  \\
\multicolumn{1}{|c}{5.3} & 10.5 &300 & 825  & 41.2  & {\bf 38.3} & 55.5  \\
\multicolumn{1}{|c}{5.3} & 10.5 &300 & 1650 & 57.5  & {\bf 53.8} & 71.2  \\
\multicolumn{1}{|c}{5.4} & 10.7 &350 & 550  & 47.2  & {\bf 44.2} & 69.0  \\
\multicolumn{1}{|c}{5.4} & 10.7 &350 & 825  & 54.3  & {\bf 48.5} & 76.3  \\
\multicolumn{1}{|c}{5.4} & 10.7 &350 & 1650 & 78.5  & {\bf 68.2} & 101.4 \\
\multicolumn{1}{|c}{5.4} & 10.8 &400 & 550  & 59.8  & {\bf 55.7} & 98.5  \\
\multicolumn{1}{|c}{5.4} & 10.8 &400 & 825  & 70.0  & {\bf 61.5} & 108.8 \\
\multicolumn{1}{|c}{5.4} & 10.8 &400 & 1650 & 105.3 & {\bf 87.2} & 142.3 \\
\hline \hline
\multicolumn{1}{|c}{1.2} & 7.4 &50  & 550  & 0.5   & 0.6  & {\bf 0.1}   \\
\multicolumn{1}{|c}{1.2} & 7.4 &50  & 825  & 0.5   & 0.5  & {\bf 0.1}   \\
\multicolumn{1}{|c}{1.2} & 7.4 &50  & 1650 & 0.5   & 0.5  & {\bf 0.1}   \\
\multicolumn{1}{|c}{1.2} & 8.0 &100 & 550  & 1.2   & 1.2  & {\bf 0.4}   \\
\multicolumn{1}{|c}{1.2} & 8.0 &100 & 825  & 1.1   & 1.1  & {\bf 0.3}   \\
\multicolumn{1}{|c}{1.2} & 8.0 &100 & 1650 & 1.1   & 1.1  & {\bf 0.3}   \\
\multicolumn{1}{|c}{1.3} & 8.5 &150 & 550  & 1.9   & 1.9  & {\bf 0.6}   \\
\multicolumn{1}{|c}{1.3} & 8.5 &150 & 825  & 1.8   & 1.8  & {\bf 0.5}   \\
\multicolumn{1}{|c}{1.3} & 8.5 &150 & 1650 & 1.7   & 1.7  & {\bf 0.5}   \\
\multicolumn{1}{|c}{1.3} & 8.7 &200 & 550  & 2.7   & 2.6  & {\bf 0.9}   \\
\multicolumn{1}{|c}{1.3} & 8.7 &200 & 825  & 2.5   & 2.5  & {\bf 0.8}   \\
\multicolumn{1}{|c}{1.3} & 8.7 &200 & 1650 & 2.4   & 2.4  & {\bf 0.7}   \\
\multicolumn{1}{|c}{1.3} & 9.0 &250 & 550  & 3.6   & 3.6  & {\bf 1.2}   \\
\multicolumn{1}{|c}{1.3} & 9.0 &250 & 825  & 3.4   & 3.4  & {\bf 1.1}   \\
\multicolumn{1}{|c}{1.3} & 9.0 &250 & 1650 & 3.2   & 3.2  & {\bf 1.0}   \\
\multicolumn{1}{|c}{1.4} & 9.1 &300 & 550  & 4.7   & 4.7  & {\bf 1.6}   \\
\multicolumn{1}{|c}{1.4} & 9.1 &300 & 825  & 4.5   & 4.5  & {\bf 1.5}   \\
\multicolumn{1}{|c}{1.4} & 9.1 &300 & 1650 & 4.1   & 4.1  & {\bf 1.3}   \\
\multicolumn{1}{|c}{1.4} & 9.3 &350 & 550  & 5.9   & 5.8  & {\bf 2.0}   \\
\multicolumn{1}{|c}{1.4} & 9.3 &350 & 825  & 5.6   & 5.6  & {\bf 1.8}   \\
\multicolumn{1}{|c}{1.4} & 9.3 &350 & 1650 & 5.0   & 5.0  & {\bf 1.6}   \\
\multicolumn{1}{|c}{1.4} & 9.4 &400 & 550  & 7.3   & 7.2  & {\bf 2.5}   \\
\multicolumn{1}{|c}{1.4} & 9.4 &400 & 825  & 6.8   & 6.8  & {\bf 2.2}   \\
\multicolumn{1}{|c}{1.4} & 9.4 &400 & 1650 & 6.2   & 6.3  & {\bf 2.0}  \\\hline
\end{tabular}
\caption{Average computing time for S-class instances.}
\label{tab:s-class_instance}
\end{table}

\begin{table}[htbp]
\centering
\footnotesize
\begin{tabular}{rlll|rrr|}
\cline{5-7}
\multicolumn{1}{l}{} &
   &
   &
  \multicolumn{1}{l|}{} &
  \multicolumn{3}{c|}{\textbf{Avg. computing time (sec)}} \\ \hline
\multicolumn{1}{|l}{\textbf{$\hat{H}$}} &
{avg max $H$} &
  \textbf{$|\tilde{A}|$} &
  \multicolumn{1}{l|}{\textbf{$u_d$}} &
  \multicolumn{1}{c}{\textbf{C-AP}} &
  \multicolumn{1}{c}{\textbf{C-UP}} &
  \multicolumn{1}{c|}{\textbf{AC-AP}} \\ \hline
\multicolumn{1}{|l}{4.7} & 9.0 & 384  & 550  & {\bf 2.9}   & 3.3   & 3.5  \\
\multicolumn{1}{|l}{4.7} & 9.0 & 384  & 825  & {\bf 3.7}   & 4.0   & 4.2  \\
\multicolumn{1}{|l}{4.7} & 9.0 & 384  & 1650 & {\bf 4.8}   & 6.5   & 5.3  \\
\multicolumn{1}{|l}{4.8} & 9.3 & 512  & 550  & {\bf 4.6}   & 4.8   & 5.7  \\
\multicolumn{1}{|l}{4.8} & 9.3 & 512  & 825  & {\bf 5.6}   & {\bf 5.6}   & 6.7  \\
\multicolumn{1}{|l}{4.8} & 9.3 & 512  & 1650 & {\bf 7.3}   & 8.7   & 8.1  \\
\multicolumn{1}{|l}{4.9} & 9.6 & 640  & 550  & {\bf 6.2}   & 6.6   & 8.1  \\
\multicolumn{1}{|l}{4.9} & 9.6 & 640  & 825  & {\bf 7.6}   & 7.7   & 9.2  \\
\multicolumn{1}{|l}{4.9} & 9.6 & 640  & 1650 & {\bf 10.4}  & 11.5  & 11.9 \\
\multicolumn{1}{|l}{5.0} & 9.8 & 800  & 550  & {\bf 9.4}   & {\bf 9.4}   & 12.2 \\
\multicolumn{1}{|l}{5.0} & 9.8 & 800  & 825  & 11.0  & {\bf 10.7}  & 14.0 \\
\multicolumn{1}{|l}{5.0} & 9.8 & 800  & 1650 & {\bf 15.1}  & 15.8  & 18.1 \\
\multicolumn{1}{|l}{5.1} & 10.0 & 959  & 550  & 12.2  & {\bf 12.1}  & 16.6 \\
\multicolumn{1}{|l}{5.1} & 10.0 & 959  & 825  & 14.3  & {\bf 13.8}  & 18.4 \\
\multicolumn{1}{|l}{5.1} & 10.0 & 959  & 1650 & {\bf 20.7}  & 21.4  & 24.8 \\
\multicolumn{1}{|l}{5.1} & 10.1 & 1151 & 550  & {\bf 15.8}  & 16.2  & 22.6 \\
\multicolumn{1}{|l}{5.1} & 10.1 & 1151 & 825  & 19.3  & {\bf 18.7}  & 25.7 \\
\multicolumn{1}{|l}{5.1} & 10.1 & 1151 & 1650 & 28.4  & {\bf 27.1}  & 35.1 \\
\multicolumn{1}{|l}{5.3} & 10.4 & 1553 & 550  & 27.4  & {\bf 26.9}  & 40.4 \\
\multicolumn{1}{|l}{5.3} & 10.4 & 1553 & 825  & {\bf 30.1}  & {\bf 30.1}  & 44.1 \\
\multicolumn{1}{|l}{5.3} & 10.4 & 1553 & 1650 & 58.0  & {\bf 54.5}  & 60.1 \\
\multicolumn{1}{|l}{5.5} & 10.6 & 1800 & 550  & 32.4  & {\bf 32.2}  & 56.9 \\
\multicolumn{1}{|l}{5.5} & 10.6 & 1800 & 825  & 364.0 &  103.3 & {\bf 65.0} \\
\multicolumn{1}{|l}{5.5} & 10.6 & 1800 & 1650 & 255.9 & 226.5 & {\bf 88.7} \\
\hline \hline
\multicolumn{1}{|l}{1.2} & 7.6 & 384  & 550  & 0.8   & 0.8   & {\bf 0.2}  \\
\multicolumn{1}{|l}{1.2} & 7.6& 384  & 825  & 0.7   & 0.7   & {\bf 0.2}  \\
\multicolumn{1}{|l}{1.2} & 7.6 & 384  & 1650 & 0.7   & 0.7   & {\bf 0.2}  \\
\multicolumn{1}{|l}{1.2} & 7.9 & 512  & 550  & 1.0   & 1.1   & {\bf 0.3}  \\
\multicolumn{1}{|l}{1.2} & 7.9 & 512  & 825  & 1.0   & 1.0   & {\bf 0.3}  \\
\multicolumn{1}{|l}{1.2} & 7.9 & 512  & 1650 & 1.0   & 1.0   & {\bf 0.2}  \\
\multicolumn{1}{|l}{1.3} & 8.2 & 640  & 550  & 1.3   & 1.3   & {\bf 0.4}  \\
\multicolumn{1}{|l}{1.3} & 8.2 & 640  & 825  & 1.2   & 1.3   & {\bf 0.3}  \\
\multicolumn{1}{|l}{1.3} & 8.2 & 640  & 1650 & 1.2   & 1.2   & {\bf 0.3}  \\
\multicolumn{1}{|l}{1.3} & 8.4 & 800  & 550  & 1.6   & 1.7   & {\bf 0.5}  \\
\multicolumn{1}{|l}{1.3} & 8.4 & 800  & 825  & 1.6   & 1.6   & {\bf 0.5}  \\
\multicolumn{1}{|l}{1.3} & 8.4 & 800  & 1650 & 1.6   & 1.6   & {\bf 0.4}  \\
\multicolumn{1}{|l}{1.3} & 8.6 & 959  & 550  & 2.0   & 2.0   & {\bf 0.7}  \\
\multicolumn{1}{|l}{1.3} & 8.6 & 959  & 825  & 1.9   & 1.9   & {\bf 0.6}  \\
\multicolumn{1}{|l}{1.3} & 8.6 & 959  & 1650 & 1.9   & 1.9   & {\bf 0.5}  \\
\multicolumn{1}{|l}{1.3} & 8.7 & 1151 & 550  & 2.5   & 2.5   & {\bf 0.9}  \\
\multicolumn{1}{|l}{1.3} & 8.7 & 1151 & 825  & 2.3   & 2.4   & {\bf 0.8}  \\
\multicolumn{1}{|l}{1.3} & 8.7 & 1151 & 1650 & 2.3   & 2.3   & {\bf 0.7}  \\
\multicolumn{1}{|l}{1.4} & 9.0 & 1553 & 550  & 3.5   & 3.6   & {\bf 1.4}  \\
\multicolumn{1}{|l}{1.4} & 9.0 & 1553 & 825  & 3.3   & 3.3   & {\bf 1.1}  \\
\multicolumn{1}{|l}{1.4} & 9.0 & 1553 & 1650 & 3.2   & 3.2   & {\bf 1.1}  \\
\multicolumn{1}{|l}{1.4} & 9.2 & 1800 & 550  & 4.3   & 4.4   & {\bf 1.9}  \\
\multicolumn{1}{|l}{1.4} & 9.2 & 1800 & 825  & 4.0   & 4.0   & {\bf 1.5}  \\
\multicolumn{1}{|l}{1.4} & 9.2 & 1800 & 1650 & 3.9   & 3.9   & {\bf 1.5}   \\ \hline
\end{tabular}
\caption{Average computing time for A-class instances}
\label{tab:a-class_instance}
\end{table}

\begin{table}[htbp]
\centering
\footnotesize
\begin{tabular}{clll|rrr|}
\cline{5-7}
\multicolumn{1}{l}{} &
   &
   &
  \multicolumn{1}{l|}{} &
  \multicolumn{3}{c|}{\textbf{Avg. computing time (sec)}} \\ \hline
\multicolumn{1}{|l}{\textbf{$\hat{H}$}} &
{avg max $H$} &
  \textbf{$|Z|$} &
  \multicolumn{1}{l|}{\textbf{$u_d$}} &
  \multicolumn{1}{c}{\textbf{C-AP}} &
  \multicolumn{1}{c}{\textbf{C-UP}} &
  \multicolumn{1}{c|}{\textbf{AC-AP}} \\ \hline
\multicolumn{1}{|c}{5.0} & 9.7 & 32 & 550  & 8.5  & {\bf 8.2}  & 11.2  \\
\multicolumn{1}{|c}{5.0} & 9.7 & 32 & 825  & 10.6 & {\bf 9.6} & 13.0  \\
\multicolumn{1}{|c}{5.0} & 9.7 & 32 & 1650 & {\bf 13.9} & 14.3 & 16.6  \\
\multicolumn{1}{|c}{5.0} & 9.7 & 48 & 550  & 9.1  & {\bf 8.4}  & 11.6 \\
\multicolumn{1}{|c}{5.0} & 9.7 & 48 & 825  & 10.5 & {\bf 9.7} & 13.2  \\
\multicolumn{1}{|c}{5.0} & 9.7 & 48 & 1650 & {\bf 14.4} & 16.4 & 16.6 \\
\multicolumn{1}{|c}{5.0} & 9.7 & 64 & 550  & {\bf 8.9}  & 9.0  & 11.5 \\
\multicolumn{1}{|c}{5.0} & 9.7 & 64 & 825  & 11.1 & {\bf 9.9} & 13.7  \\
\multicolumn{1}{|c}{5.0} & 9.7 & 64 & 1650 & {\bf 14.3} & 17.8 & 16.7 \\
\multicolumn{1}{|c}{5.0} & 9.7 & 96 & 550  & {\bf 9.0}  & 9.8 & 11.6  \\
\multicolumn{1}{|c}{5.0} & 9.7 & 96 & 825  & {\bf 11.1} & 11.2 & 13.4\\
\multicolumn{1}{|c}{5.0} & 9.7 & 96 & 1650 & {\bf 14.1} & 23.5 & 16.5  \\
\hline \hline
\multicolumn{1}{|c}{1.3} & 8.3 & 32 & 550  & 1.6  & 1.6  & {\bf 0.5}   \\
\multicolumn{1}{|c}{1.3} & 8.3 & 32 & 825  & 1.5  & 1.5  & {\bf 0.4}   \\
\multicolumn{1}{|c}{1.3} & 8.3 & 32 & 1650 & 1.5  & 1.5  & {\bf 0.4}   \\
\multicolumn{1}{|c}{1.3} & 8.3 & 48 & 550  & 1.6  & 1.7  & {\bf 0.5}   \\
\multicolumn{1}{|c}{1.3} & 8.3 & 48 & 825  & 1.5  & 1.6  & {\bf 0.4}   \\
\multicolumn{1}{|c}{1.3} & 8.3 & 48 & 1650 & 1.6  & 1.5  & {\bf 0.4}   \\
\multicolumn{1}{|c}{1.3} & 8.3 & 64 & 550  & 1.5  & 1.7  & {\bf 0.5}   \\
\multicolumn{1}{|c}{1.3} & 8.3 & 64 & 825  & 1.5  & 1.6  & {\bf 0.4}   \\
\multicolumn{1}{|c}{1.3} & 8.3 & 64 & 1650 & 1.5  & 1.5  & {\bf 0.4}   \\
\multicolumn{1}{|c}{1.3} & 8.3 & 96 & 550  & 1.6  & 1.9  & {\bf 0.5}   \\
\multicolumn{1}{|c}{1.3} & 8.3 & 96 & 825  & 1.6  & 1.7  & {\bf 0.4}   \\
\multicolumn{1}{|c}{1.3} & 8.3 & 96 & 1650 & 1.6  & 1.7  & {\bf 0.4}    \\ \hline
\end{tabular}
\caption{Average computing time for Z-class instances}
\label{tab:z-class_instance}
\end{table}

\begin{table}[htbp]
\centering
\footnotesize
\begin{tabular}{rlll|rrr|}
\cline{5-7}
\multicolumn{1}{l}{} &
   &
   &
  \multicolumn{1}{l|}{} &
  \multicolumn{3}{c|}{\textbf{Avg. computing time (sec)}} \\ \hline
\multicolumn{1}{|l}{\textbf{$\hat{H}$}} &
{avg max $H$} &
  \textbf{$|D|$} &
  \multicolumn{1}{l|}{\textbf{$u_d$}} &
  \multicolumn{1}{c}{\textbf{C-AP}} &
  \multicolumn{1}{c}{\textbf{C-UP}} &
  \multicolumn{1}{c|}{\textbf{AC-AP}}\\ \hline
\multicolumn{1}{|c}{5.4} & 9.5 & 21 & 393  & {\bf 10.2}  & 10.8   & 13.0     \\
\multicolumn{1}{|c}{5.4} & 9.5 & 21 & 589  & {\bf 13.7}  & 13.9   & 16.7    \\
\multicolumn{1}{|c}{5.4} & 9.5 & 21 & 1179 & {\bf 14.3}  & 24.4   & 17.4    \\
\multicolumn{1}{|c}{5.9} & 9.5 & 27 & 306  & {\bf 13.1}  & 14.7   & 17.9    \\
\multicolumn{1}{|c}{5.9} & 9.5 & 27 & 458  & {\bf 16.3}  & 19.2   & 20.3    \\
\multicolumn{1}{|c}{5.9} & 9.5 & 27 & 917  & {\bf 43.6}  & 131.1  & 54.7   \\
\multicolumn{1}{|c}{6.2} & 9.6 & 33 & 250  & {\bf 15.2}  & 21.4   & 22.3    \\
\multicolumn{1}{|c}{6.2} & 9.6 & 33 & 375  & {\bf 55.5}  & 63.2   & 58.7    \\
\multicolumn{1}{|c}{6.2} & 9.6 & 33 & 750  & {\bf 288.1} & $1141.1^{(1)}$ & 294.4 \\
\multicolumn{1}{|c}{6.5} & 9.6 & 39 & 212  & {\bf 20.8}  & 30.1   & 31.4    \\
\multicolumn{1}{|c}{6.5} & 9.6 & 39 & 317  & {\bf 46.7}  & 64.4   & 54.6   \\
\multicolumn{1}{|c}{6.5} & 9.6 & 39 & 635  & 211.4 & $1240.4^{(2)}$ & {\bf 205.0}       \\
\multicolumn{1}{|c}{6.7} & 9.7 & 45 & 183  & {\bf 27.1}  & 40.6   & 41.6   \\
\multicolumn{1}{|c}{6.7} & 9.7 & 45 & 275  &  70.5  & 98.5  & {\bf 69.3}   \\
\multicolumn{1}{|c}{6.7} & 9.7 & 45 & 550  & 925.9 & ${\bf 669.0}^{(3)}$  & 900.8      \\
\hline \hline
\multicolumn{1}{|c}{1.7} & 8.3 & 21 & 393  & 1.7   & 1.7    & {\bf 0.6}       \\
\multicolumn{1}{|c}{1.7} & 8.3 & 21 & 589  & 1.6   & 1.7    & {\bf 0.6}       \\
\multicolumn{1}{|c}{1.7} & 8.3 & 21 & 1179 & 1.6   & 1.6    & {\bf 0.5}       \\
\multicolumn{1}{|c}{2.0} & 8.4 & 27 & 306  & 2.5   & 2.5    & {\bf 1.1}       \\
\multicolumn{1}{|c}{2.0} & 8.4 & 27 & 458  & 2.5   & 2.5    & {\bf 1.0}       \\
\multicolumn{1}{|c}{2.0} & 8.4 & 27 & 917  & 2.4   & 2.3    & {\bf 0.9}       \\
\multicolumn{1}{|c}{2.3} & 8.6 & 33 & 250  & 3.5   & 3.6    & {\bf 1.6}       \\
\multicolumn{1}{|c}{2.3} & 8.6 & 33 & 375  & 3.6   & 3.6    & {\bf 1.6}       \\
\multicolumn{1}{|c}{2.3} & 8.6 & 33 & 750  & 3.2   & 3.2    & {\bf 1.3}       \\
\multicolumn{1}{|c}{2.6} & 8.7 & 39 & 212  & 4.5   & 4.5    & {\bf 2.2}       \\
\multicolumn{1}{|c}{2.6} & 8.7 & 39 & 317  & 4.7   & 4.8    & {\bf 2.3}       \\
\multicolumn{1}{|c}{2.6} & 8.7 & 39 & 635  & 4.4   & 4.4    & {\bf 1.9}       \\
\multicolumn{1}{|c}{2.8} & 8.8 & 45 & 183  & 5.7   & 5.7    & {\bf 2.8}       \\
\multicolumn{1}{|c}{2.8} & 8.8 & 45 & 275  & 5.8   & 5.9    & {\bf 2.9}       \\
\multicolumn{1}{|c}{2.8} & 8.8 & 45 & 550  & 5.9   & 6.0    & {\bf 2.7}       \\ \hline
\multicolumn{7}{l}{(1) 9/10 solved to optimality, 0.04\% gap for 1/10} \\
\multicolumn{7}{l}{(2) 9/10 solved to optimality, 0.02\% gap for 1/10} \\
\multicolumn{7}{l}{(3) 1/10 solved to optimality, 0.06\% average gap for 9/10} \\
\end{tabular}
\caption{Average computing time for D-class instances}
\label{tab:d-class_instance}
\end{table}

\section{Case Study}\label{sec:study_case}

This section presents a case study where the LP seeks to implement CDPs in the region of Montreal, Quebec. For this study, we utilize data from our industrial partner, Purolator. We transform data related to cost, demand, and preferences to preserve confidentiality. In this case, the demand is known at the postal codes, which are grouped in $35$ FSA (forward sorting area). Thus, $|Z|=35$, considering only those FSAs that are served from the main terminal. In addition, only postal codes having regular (non-zero) weekly demand are considered leading to $|A|=1,081$ postal codes. The total weekly demand amounts to 29,092 parcels. We assume three categories of customers, and we make them fictitious for confidentiality reasons. Table~\ref{tab:case_categories} shows the three categories of customers and the corresponding utility parameters. The utility only depends on the distance from CDPs to the customer's corresponding location ($\alpha^2_{k}=0$).  In addition, we assume that the utility function represents an accurate estimation of customer preferences with $\beta=0.25$.
Nevertheless, we illustrate the impact of using $\beta=1$ which implies a higher degree of uncertainty. We also analyze the effect of more sensitive customers. For this case study, we impose a 10-hour computing time budget. 

\begin{table}[htbp]
\centering
\small
\begin{tabular}{lcc}
\hline
\textbf{Customer category} & \textbf{\% Population} & \textbf{$\alpha_k^1$}  \\ \hline
Category 1 & 20\% & -0.10  \\
Category 2     & 70\% & -0.15    \\
Category 3                          & 10\% & -0.25    \\ \hline
\end{tabular}
\caption{Customer categories and utility function parameters}
\label{tab:case_categories}
\end{table}

\begin{figure}[htbp]
    \centering
    \includegraphics[width=1\textwidth]   {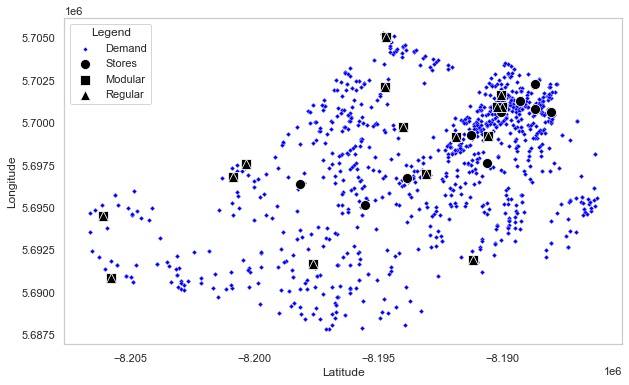}
    \caption{Postal Codes and CDPs.}
    \label{fig:pcodes_cdp}
\end{figure}

We consider 40 CDP candidate locations, where 10 correspond to stores ($D^\text{P}$), 15 to modular parcel lockers, and 15 to regular parcel lockers $(D^\text{R})$. In Table~\ref{tab:CDP_case_study}, we present the data related to the CDPs. Note that CDPs at stores charge an additional $\$1$ for delivering parcels to customers. In addition, $f_d$ correspond to weekly amortized fixed costs, and $u_d$ corresponds to the capacity of CDPs for a week. Figure~\ref{fig:pcodes_cdp} shows the centroid location of postal codes and the potential CDPs for this case. 

Without considering CDPs, the (current) estimated total routing cost per week is \$73,301. This cost is reduced by 3.1\% in the optimal solution. In this case, 7 regular parcel lockers are open, capturing 26.4\% of the demand. Noteworthy is that even if modular parcel lockers display a better unit cost capacity ($f_d/u_d$) than regular parcel lockers, they are not selected. However, the much lower capacity of modular parcel lockers limits the total cost reductions they can provide.

If, instead, customers in the first and second categories are 25\% more sensitive to distance, only 6 regular parcel lockers are needed, reducing the current estimated total costs by 2.3\% and capturing 21.9\% of the demand. On the contrary, a choice model with $\beta=1$ significantly overestimates the demand captured by CDPs (56\%) which leads to opening 14 CDPs.

We note that regular parcel lockers are preferred over the other two classes of CDPs, but this is mainly due to their associated costs and capacity. If the cost for delivering a parcel by a store were $\$0.5$ instead of $\$1$, then 9 CDPs are open at stores, and 5 regular parcel lockers are implemented, resulting in a routing cost reduction of 3.7\%. In addition, in such a case, the 14 CDPs capture 36.1\% of the demand.

In summary, CDPs' capacities not only limit the demand they can attract but also negatively impact the potential cost reductions they offer. In addition, underestimating customers' negative perception of distance (leading to higher uncertainty) may lead to an overestimation of the number of CDPs. 

\begin{table}[htbp]
\centering
\small
\begin{tabular}{lrrrr}
\hline
\textbf{Type of CDP} &
  \multicolumn{1}{c}{\textbf{$u_d$}} &
  \multicolumn{1}{c}{\textbf{$f_d$}} &
  \multicolumn{1}{c}{\textbf{$b_d$}} &
  \multicolumn{1}{c}{\textbf{$c_{dz}^\text{AR}$}} \\ \hline
Stores  & 600.0 & 0.0   & 300.0 & $0.7 \times c_z^\text{R} +1$ \\
Modular & 148 & 31.0  & 0.0   & $0.7 \times c_z^\text{R}$  \\
Regular   & 1200  & 500.0 & 0.0   & $0.7 \times c_z^\text{R}$   \\ \hline
\end{tabular}
\caption{CDP instance data}
\label{tab:CDP_case_study}
\end{table}

\section{Conclusions}\label{sec:conclusion}
We addressed a choice-based capacitated facility location problem with random utility-maximizing customers. We considered the case of an LP that seeks to locate CDPs of different types. At this strategic level, the LP has imperfect knowledge of customer preferences as well as precise customer locations. The CDP locations and, hence, the benefit for the LP of these facilities depend on the capacity to attract customers. Thus, customer locations and preferences for the CDPs are key drivers for the decisions. 

We proposed two formulations including closest assignment constraints that we solved using a standard Benders decomposition method. We showed that the closest assignment constraints are derived from a single-level reformulation of a bilevel program with a specific formulation for the follower's problem.
We conducted an extensive experimental study to analyze the impact of the level of uncertainty (customer preferences and geographical distribution) on the optimal solutions and computing time. We used an entropy measure for interpretability. Consistent with the finding in \cite{legault2023modelfree}, the most \rev{efficient} model and solution approach depend on the level of uncertainty in the instance. This is captured by the entropy and, importantly, it can be easily computed before solving an instance. Using scenario aggregation is the best option for low-entropy instances. For higher entropy levels, the formulation without aggregation resulted in shorter computing times. In  this case, we considered two ways to define the subproblems: the standard (referred to as automatic partition) way where integer variables are in the master and continuous variables in the subproblems, and user partition where we also include continous variables in the master resulting in exclusively generating feasibility cuts. The user partition worked particularly well when increasing the number of scenarios, whereas the automatic partition displayed strong performance more consistently across all instances. 

The experimental results showed that the results are quite sensitive to the demand model and accuracy of the geographical distribution of demand. This highlights the need for high-performing demand models, especially in regions where it impacts the decisions (and objective function) the most. \rev{The proposed methodology can accommodate any random utility maximization model whose utility is linear in the decision-variables. Future work should focus on estimating such models and empirically investigating the impact on the resulting decisions. In this context,} we believe that integrated learning and optimization \citep{SADANA2025271} for training such models constitutes a \rev{challenging yet} potentially  valuable direction for future research.

In this work we assumed that demand is allocated to home delivery if there is insufficient capacity at the preferred CDP. Relaxing this assumption would lead to a more complex setting which we intend to explore in future research.

\section*{Author Contribution}
\author{\bf{David Pinzon}}: Conceptualization, Methodology, Software, Validation, Investigation, Formal analysis, Writing – Original Draft. \author{\bf{Ammar Metnani}}: Conceptualization, Supervision, Review \& Editing. \author{\bf{Emma Frejinger}}: Supervision, Conceptualization, Resources, Validation, Writing – Review \& Editing, Funding acquisition. 

\section*{Acknowledgements}
We gratefully acknowledge that the first author received financial support from the Ecuadorian Secretar\'ia de Educaci\'on Superior, Ciencia, Tecnolog\'ia e Innovacion (SENESCYT) under its international scholarships program and from University Escuela Superior Polit\'ecnica del Litoral, Ecuador. The research is funded by Data Intelligence for Logistics Research Chair at Universit\'e de Montr\'eal thanks to support from Purolator and Natural Sciences and Engineering Research Council of Canada (grant CRDPJ 538506 – 19). We express our deep gratitude to Eric Larsen who helped us improve the paper. The first author got financial support from EURO to attend the EURO Summer Institute in Location Science, held in Edinburgh in June 11-22, 2022. Ideas from that event led to the formulation with closest assignment constraints that we propose in this work.

\bibliography{article2}

\appendix

\clearpage
    \section{Notation}\label{art2:Main_notation}
        \begin{table}[h!]
        \small
        \centering
        \begin{tabular}{|ll|}
        \hline
        \multicolumn{2}{|l|}{\textbf{Indexes}} \\ 
        $d$            & CDP\\ 
        $z$            & Zone \\ 
        $a$            & Subzone \\
        $k$            & Customer category \\
        $q$            & pattern index \\ 
        $l$            & capacity level \\
        &\\
        \multicolumn{2}{|l|}{\textbf{Sets}} \\ 
        $A$            & Set of Subzones \\
        $A_z$          & Set of subzones in zone $z$ \\
        $Z$            & Set of Zones  \\ 
        $K$            & Set of Customer categories \\
        $L$            & Set of Capacity levels  \\
        $D$            & Set of CDPs   \\ 
        $D^\text{P}$        & Set of pick up points or stores \\
        $D^\text{M}$        &Set of  modular parcel lockers \\ 
        $D^\text{R}$        & Set of regular parcel lockers\\ 
        $Q$             & Set of Patterns \\
        & \\
        \multicolumn{2}{|l|}{\textbf{Deterministic parameters}} \\ $c_{z}^\text{R}$    & Estimated routing cost for zone $z$ \\ 
        $c_{dz}^{\text{AR}}$    & Estimated routing cost for zone $z$ affected by CDP $d$\\ 
        $u_{d}$         & Capacity of CDP $d$   \\ 
        $f_{d}$         & Fixed cost of CDP $d$  \\ 
        $u_d^l$         & Capacity level $l$ of CDP $d$   \\ 
        $f_d^l$         & Fixed cost associated to capacity level $l$ of CDP $d$  \\ 
        $\underline{b}_d$       & Minimum demand level required for CDP $d \in D^\text{P}$ \\
        $b_{z}$        & Average number of parcels to deliver in zone $z$    \\  
        $b_{az}$    & Average number of parcels to deliver in subzone $a$ \\ 
        $\theta_a$    & Coordinates of the centroid of subzone $a$ \\
        $\theta_d$    & Coordinates of the location of CDP $d$ \\
        $b_{zak}$ & average number of parcels of customers in category $k$ from subzone $a$ in zone $z$ \\
        $\alpha_k^1$ & Distance-related coefficient in the utility function\\
        $\alpha_k^2$ & Coefficient in the utility function related to exogenous characteristics\\
        $e_{gdask}$ & 1, If $U_{gaks} \geq U_{daks}$; 0, otherwise\\
        \hline
        \end{tabular}
        \caption{Indexes, sets and deterministic parameters}
        \label{tab:art2:notation_part1}
        \end{table}
        
        \clearpage

        \begin{table}[h!]
        \small
        \centering
        \begin{tabular}{|ll|}
        \hline
        \multicolumn{2}{|l|}{\textbf{Derived variables}} \\ 
        $\rho_{dak}(\boldsymbol{x})$ &  Probability that CDP $d$ captures demand from customers in category $k$ and subzone $a$  \\
        $\psi_{dz}(\boldsymbol{x})$ & Probability that CDP $d$ captures demand from zone $z$\\
        $\hat{\rho}_{dak}(\boldsymbol{x})$ &SAA of $\psi_{dak}(\boldsymbol{x})$ \\
        $\nu_{qak}$ & Fraction of customers in subzone $a$ and category $k$ with the same preference pattern $q$ \\
        $U_{dak}$ & Random utility for CDP $d$ an customers in category $k$ and subzone $a$ \\
        $U_{daks}$ & Realization of $U_{dak}$ for scenario $s$\\
        $\varepsilon_{dak}$ & Random term in the utility function\\
        $\varepsilon_{s}$ & Realization of $\varepsilon_{dak}$ in scenario $s$\\
        \hline
        \multicolumn{2}{|l|}{\textbf{Decision variables}} \\ $x_d$    & 1, If CDP $d$ is open; 0, otherwise \\ 
        $w_{daks}$ & 1, If customers in category $k$ and subzone $a$ are assigned to CDP $d$ in scenario $s$; 0, otherwise \\
        $p_{dz}$ & Effective fraction of demand in zone $z$ captured by CDP $d$ given the capacity restrictions \\
        $p_{0z}$ & Effective fraction of demand in zone $z$ captured by or assigned to home delivery service \\
        $w_{qd}$ & 1; If demand in pattern $q$ is assigned to CDP $d$; 0, otherwise \\
        $r_d^l$ & 1, If capacity level $ld$ is assigned to CDP $d$; 0, otherwise
         \\ 
        \hline
        \end{tabular}
        \caption{Derived and decision variables}
        \label{tab:art2:notation_part2}
        \end{table}

   \clearpage
    \section{Dominated Single-level Reformulation}\label{apx:DSLR}
    This section introduces the MILP model resulting from reformulating \eqref{af:P1} by using linear programming strong duality, a similar reasoning as in Proposition~\ref{CAC-derive}, and standard linearization techniques: 
        \begin{subequations}\label{af:RFP1}
        \begin{equation}
            \text{min } \sum_{z \in Z}  c_{z}^\text{R}b_zp_{0z}+  \sum_{z \in Z}\sum_{d \in D} c_{dz}^\text{AR}b_zp_{dz} + \sum_{d \in D^\text{M}} \sum_{l \in L_d} f^l_d r^l_d +\sum_{d \in D^\text{A}} f_d x_d  \label{SL1:of} 
        \end{equation}
         \begin{align}
             \text{s.t: } & (\ref{enlm:1}- \ref{enlm:2.8}, \ref{enlm:4} - \ref{enlm:dom_r}) \nonumber \\
             & p_{dz} \leq \frac{1}{|S|}\sum_{s \in S}\sum_{a \in A_{z}}\sum_{k \in K}\frac{b_{zak}}{b_z} v_{daks}, &  d \in D, z \in Z \label{eq:FP1.bound_prob} \\
             & p_{0z} \geq \frac{1}{|S|}\sum_{s \in S}\sum_{a \in A_{z}}\sum_{k \in K}\frac{b_{zak}}{b_z} v_{0aks}, &  z \in Z \label{eq:FP1.bound_prob_0} \\
           & \sum_{g \in D} U_{gaks} v_{gaks} +U_{0aks}w_{0aks} \geq U_{daks} x_d, & d \in D, a\in A, k \in K, s \in S \label{eq:SLFP1.1} \\
           &\sum_{g \in D} U_{gaks} v_{gaks} +U_{0aks}w_{0aks} \geq U_{0aks},  & a\in A, k \in K, s \in S \label{eq:SLFP1.2}\\
            &\sum_{d \in D \cup \{0\}} w_{daks} = 1, & a \in A_, k\in K, s\in S \label{eq:FP1.assign} \\
            &v_{daks} \leq x_d, & d \in D, a\in A, k \in K, s \in S \label{eq:F1.lr1} \\
            &v_{daks} \leq w_{daks}, &d \in D, a\in A, k \in K, s \in S \label{eq:F1.lr2}\\
            &v_{daks} \geq x_d +w_{daks}-1, & d \in D, a \in A, k \in K, s \in S \label{eq:F1,lr3} \\
            &v_{daks} \geq 0, & d \in D, a\in A, k \in K, s \in S \label{eq:F1.lr4}
        \end{align}
        \end{subequations}

    Note that this formulation differs from model formulation~\eqref{model:Sim_CAC} in Constraints~\eqref{eq:SLFP1.1}, \eqref{eq:SLFP1.2}, \eqref{eq:F1.lr1} - \eqref{eq:F1,lr3}. Constraints~\eqref{eq:SLFP1.1} and \eqref{eq:SLFP1.2} force the assignment to the delivery service with the highest utility. Constraints~\eqref{eq:F1.lr1} - \eqref{eq:F1,lr3} allow to linearize $v_{daks}=w_{adks}x_d$. 

    We use two approaches to solve model~\eqref{af:RFP1}, the fist, denoted by D-BB, applies branch-and-bound whereas the second, denoted by D-B, applies Benders decomposition.

\section{Additional Results}\label{apx:Add_results}
\begin{figure}[htbp]
    \centering
    \subfloat[Medium entropy\label{fig:cdps_frac_med_entropy}]{{\includegraphics[width=0.55\textwidth]{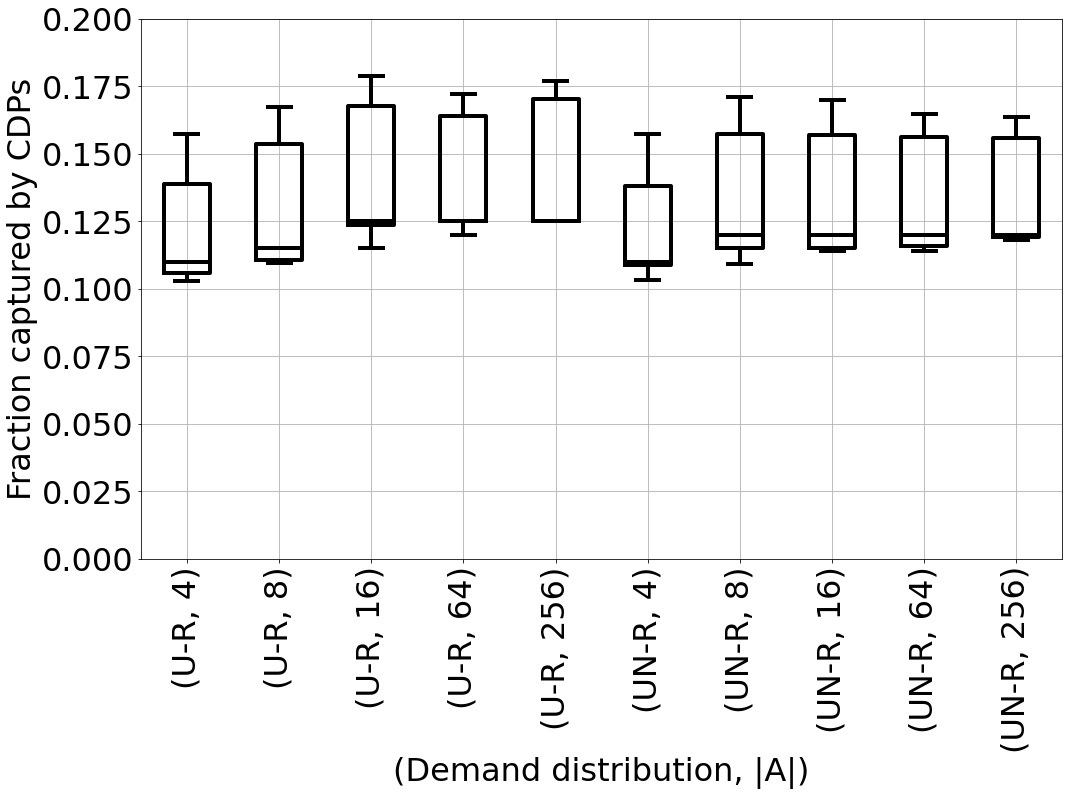}}}
    \subfloat[Low entropy \label{fig:cdps_frac_low_entropy} ]{{\includegraphics[width=0.55\linewidth]{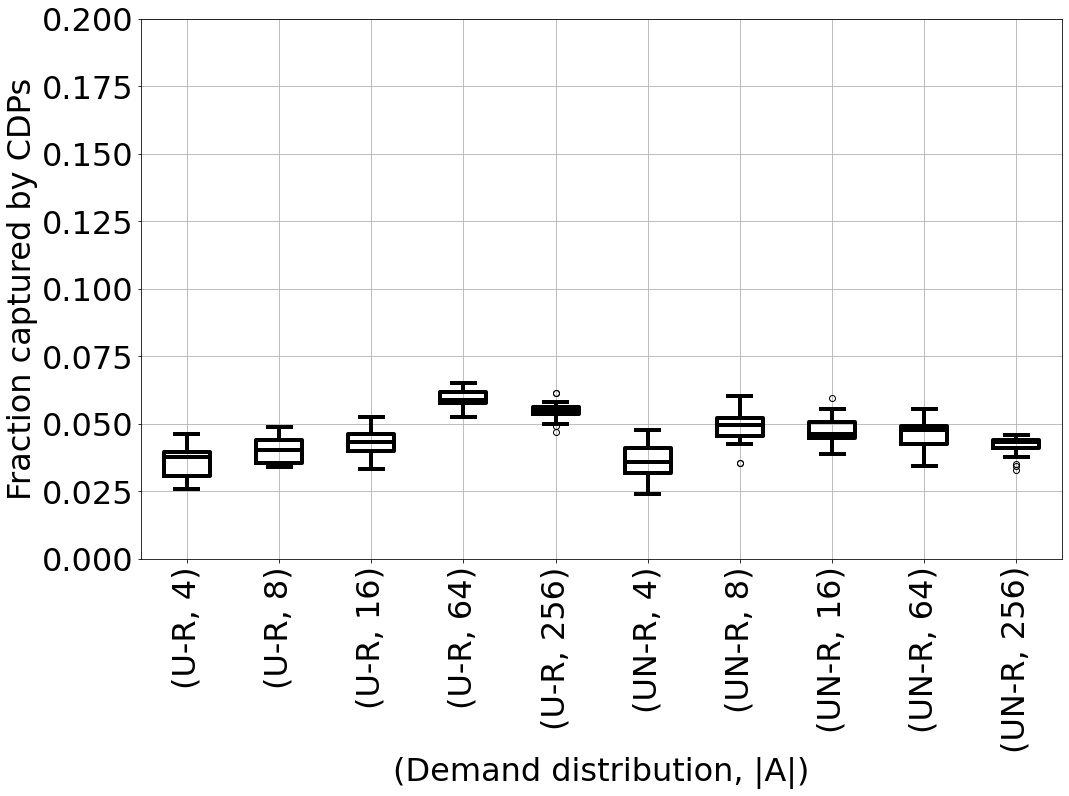}}}
    \caption{Distribution of optimal fraction captured by CDPs}
    \label{fig:Cdps_frac_cap}
\end{figure}

\begin{figure}[htbp]
    \centering
    \subfloat[C-AP vs AC-AP\label{fig:ratio_time_C-AP_vs_AC-Ap}]{{\includegraphics[width=0.55\textwidth]{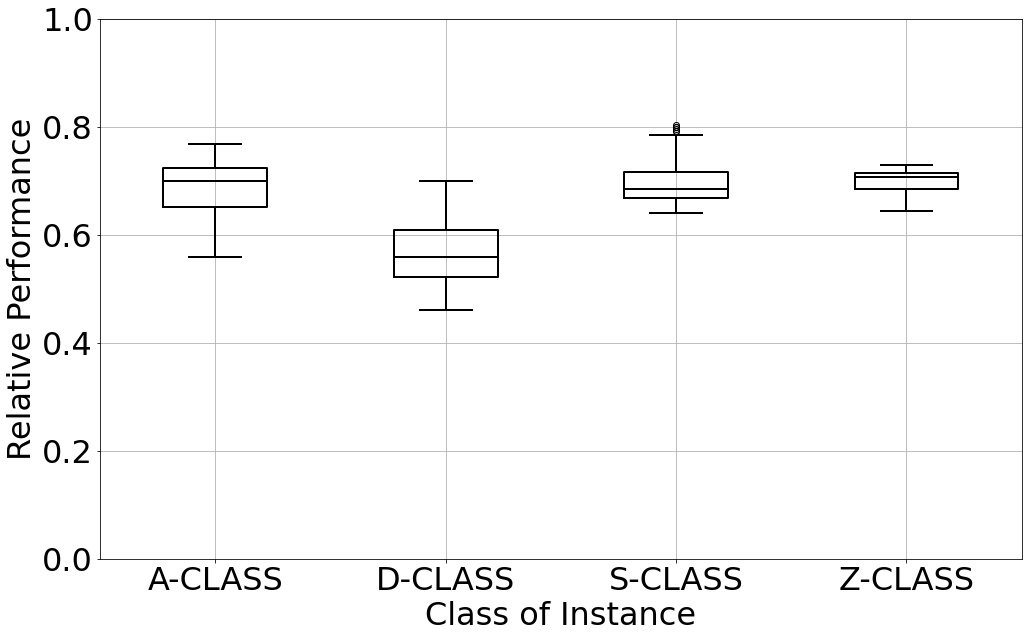}}}
    \subfloat[C-UP vs AC-AP \label{fig:ratio_time_C-UP_vs_AC-Ap} ]{{\includegraphics[width=0.55\linewidth]{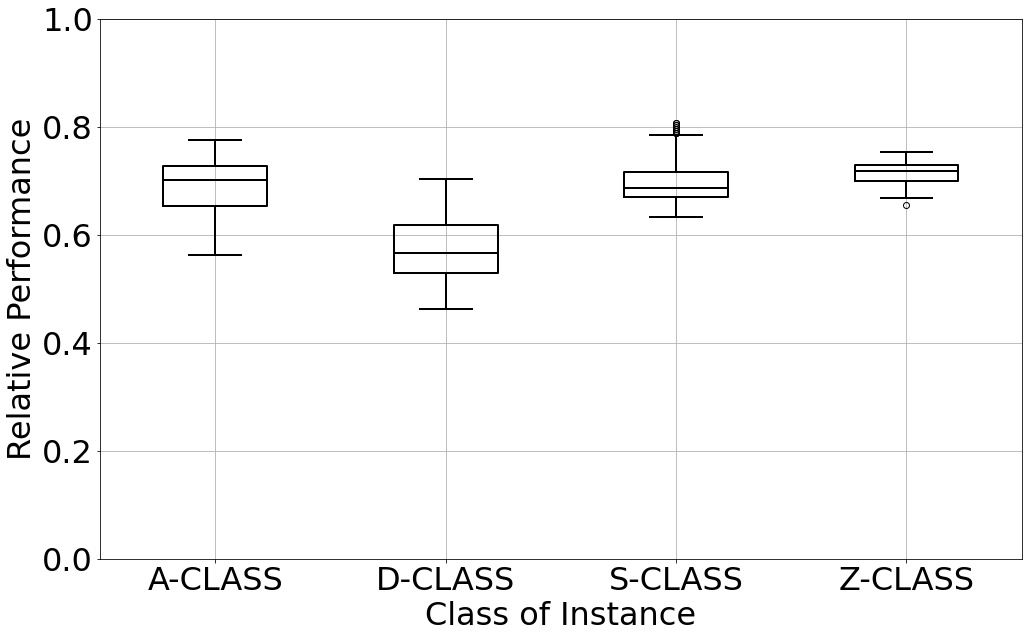}}}
    \caption{Distribution of the relative performance in low entropy instances}
    \label{fig:ratio_time_ac_ap}
\end{figure}

This section provides additional results about performance and optimal solutions.
Table~\ref{tab:com_baseline} reports the computing times for D-B, D-BB and C-AP, for a subset of instances used in \ref{sec:Interpret}. The panels of Figure~\ref{fig:Cdps_frac_cap} show the distribution of the effective fraction captured by CDPs. \rev{The} panels of Figure~\ref{fig:ratio_time_ac_ap} show the distribution of the relative performance of C-AP and C-UP with respect to AC-AP. 

\begin{figure}[htbp]
    \centering
    \subfloat[Medium entropy \label{fig:cost_M_550}]{{\includegraphics[width=0.55\textwidth]{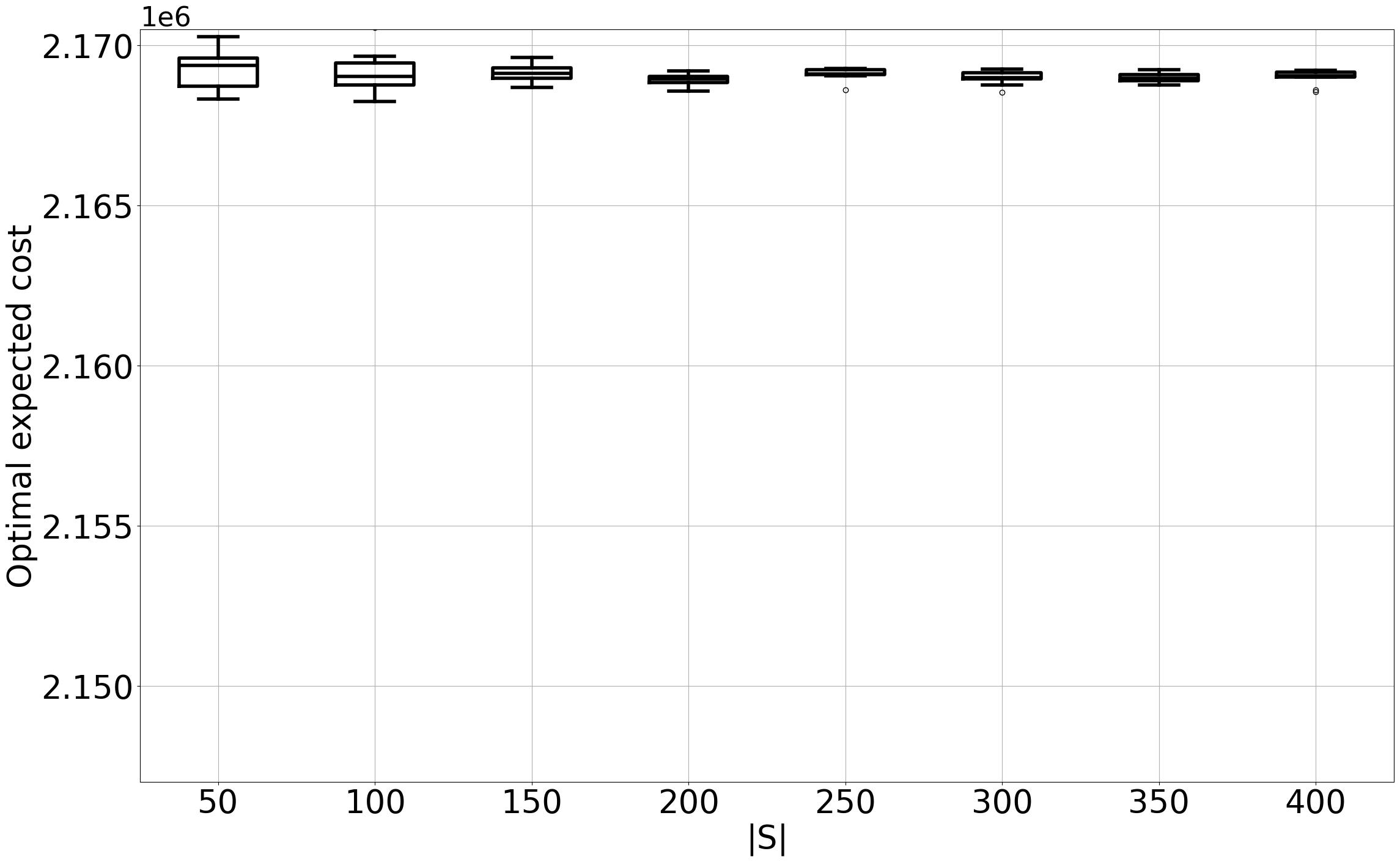}}}
    \subfloat[Low entropy \label{fig:cost_L_550} ]{{\includegraphics[width=0.55\linewidth]{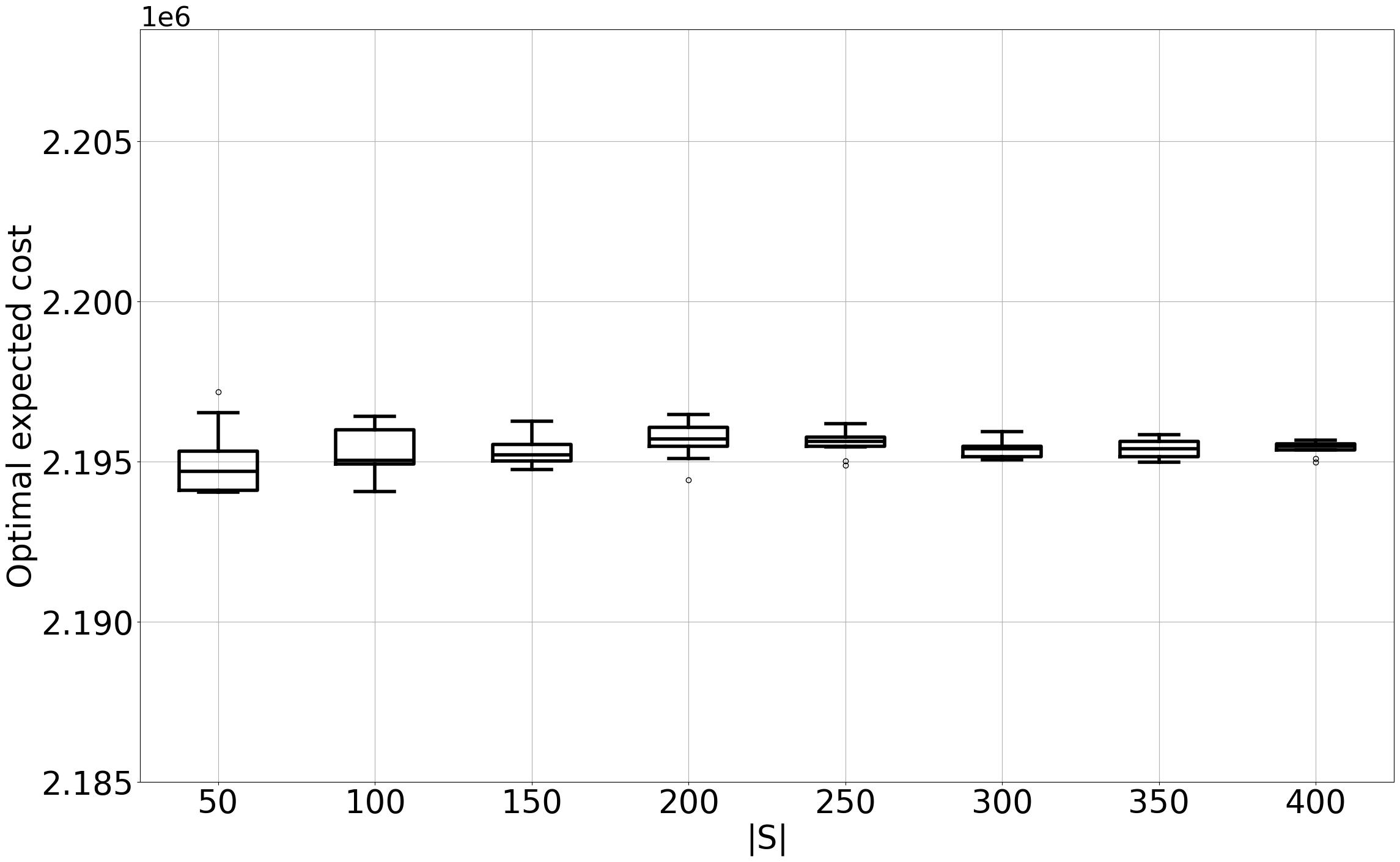}}}
    \caption{Distribution of optimal expected cost for $u_d=550$ and $|S|$}
    \label{fig:costs_550}
\end{figure}

\begin{figure}[htbp]
    \centering
    \subfloat[Medium entropy \label{fig:cost_M_825}]{{\includegraphics[width=0.55\textwidth]{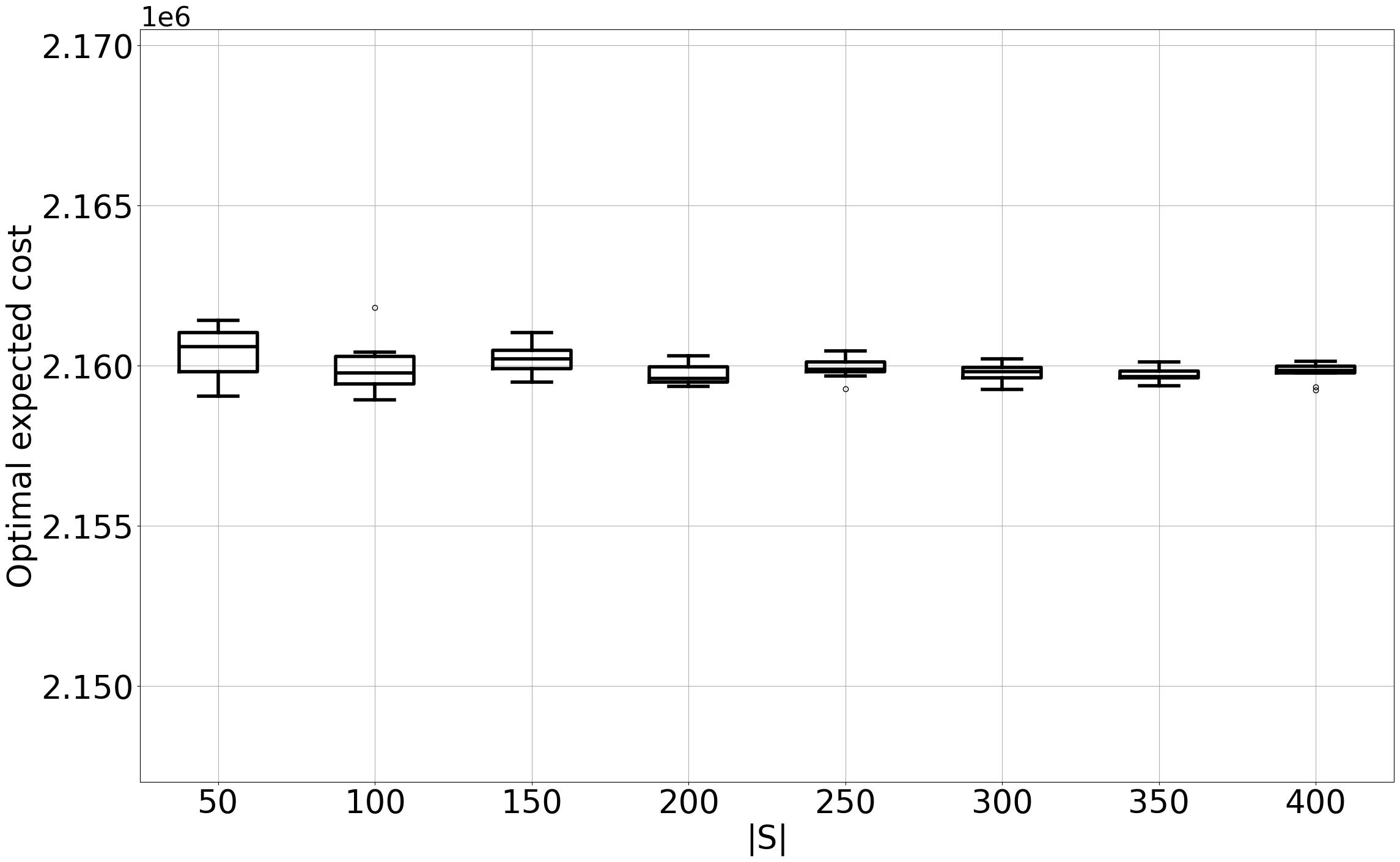}}}
    \subfloat[Low entropy \label{fig:cost_L_825} ]{{\includegraphics[width=0.55\linewidth]{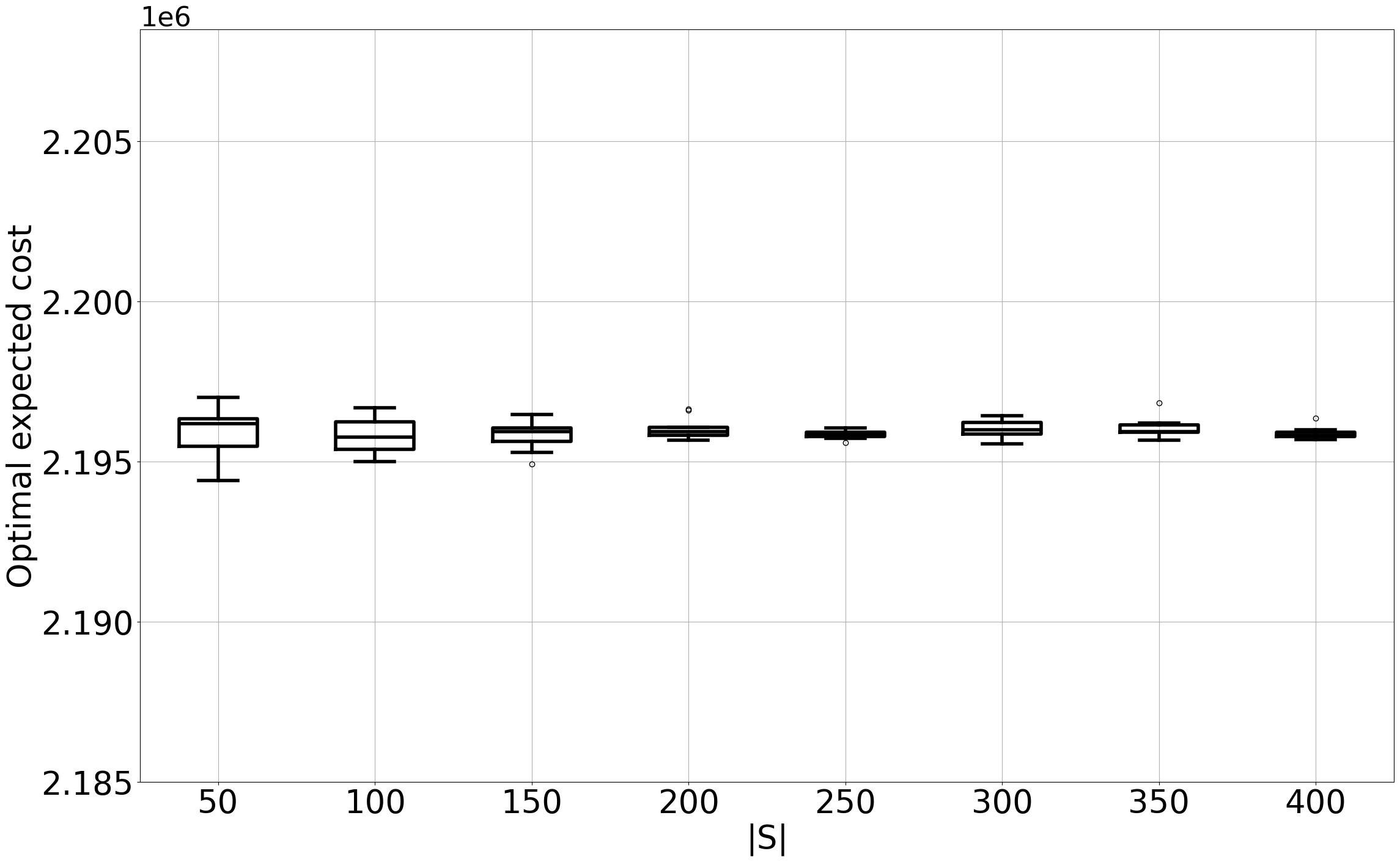}}}
    \caption{Distribution of optimal expected cost for $u_d=825$ and $|S|$}
    \label{fig:costs_825}
\end{figure}

\begin{figure}[htbp]
    \centering
    \subfloat[Medium entropy \label{fig:cost_M_1650}]{{\includegraphics[width=0.55\textwidth]{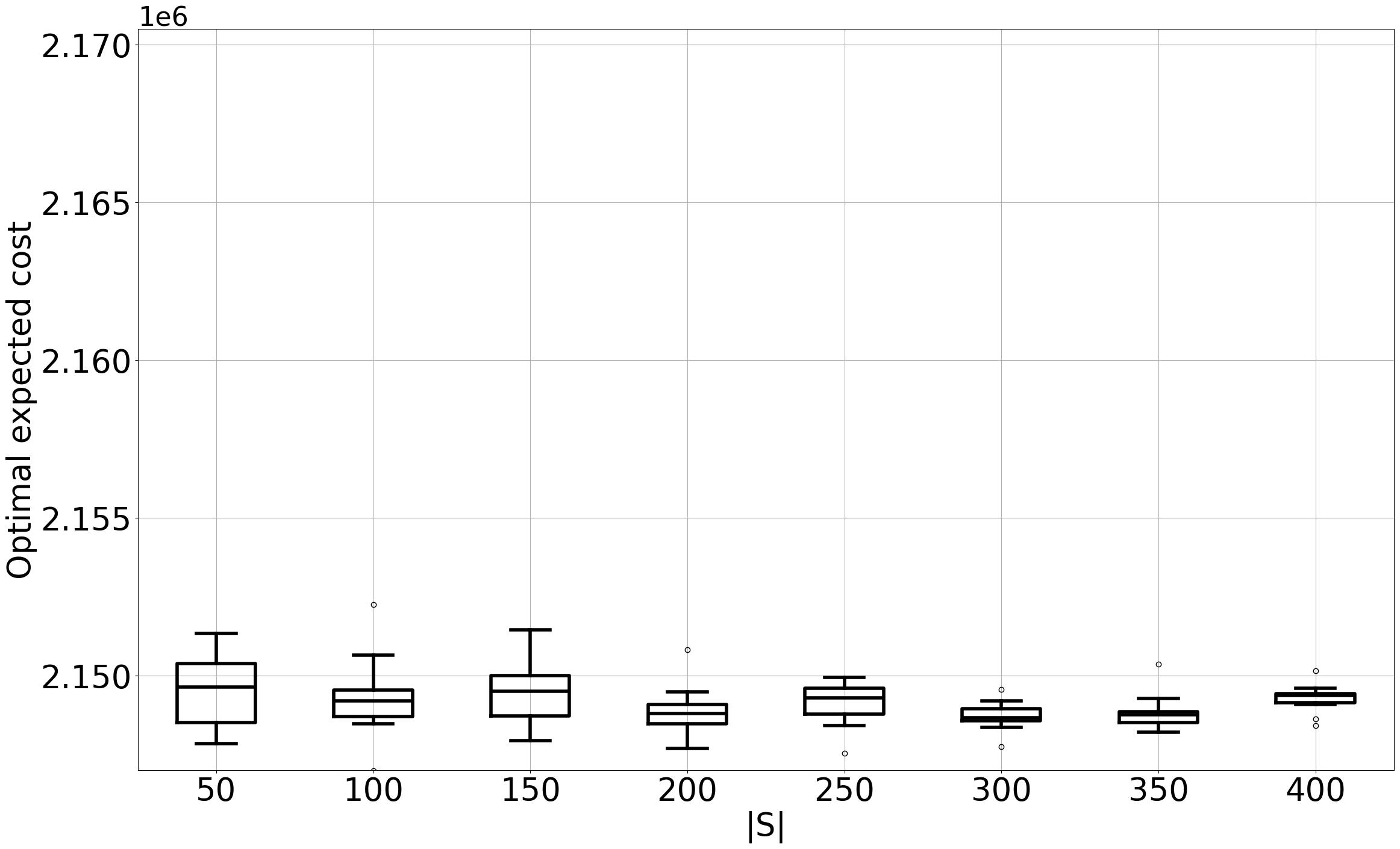}}}
    \subfloat[Low entropy \label{fig:cost_L_1650} ]{{\includegraphics[width=0.55\linewidth]{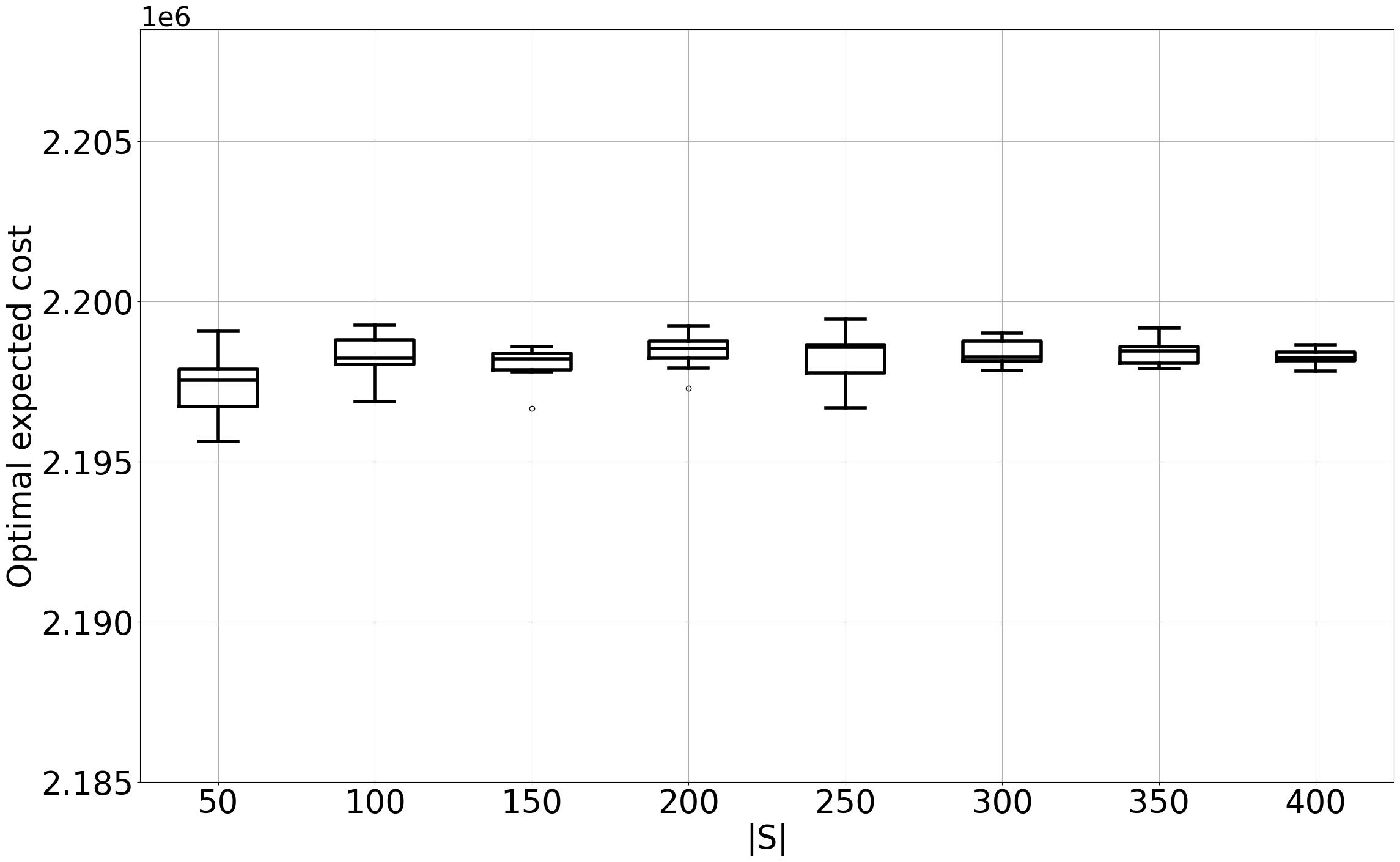}}}
    \caption{Distribution of optimal expected cost for $u_d=1650$ and $|S|$}
    \label{fig:costs_1650}
\end{figure}

\rev{The panels of Figures~\ref{fig:costs_550}, ~\ref{fig:costs_825} and~\ref{fig:costs_1650} show the distribution of the optimal expected cost for increasing $|S|$ for the three different capacity levels, respectively. As expected, increasing the size of the set of scenarios leads to a more accurate estimation of the objective function values but we note that the distributions for 50 scenarios are relatively close to those for higher number of scenarios.} 


\begin{table}[htpb]
\centering
\footnotesize
\begin{tabular}{|lllr|rrr|}
\hline
\textbf{Distribution} &
  \textbf{$|A|$} &
  \textbf{$u_d$} &
  \multicolumn{1}{l|}{\textbf{Avg. Entropy}} &
  \multicolumn{1}{l}{\textbf{D-BB}} &
  \multicolumn{1}{l}{\textbf{D-B}} &
  \multicolumn{1}{l|}{\textbf{C-AP}} \\ \hline
U-R & 4  & 20 & 5.6 & 1.1   & 1.1    & 0.1 \\
U-R & 4  & 20 & 2.9 & 2.5   & 5.9    & 0.0 \\
U-R & 4  & 20 & 0.7 & 0.5   & 0.8    & 0.0 \\
U-R & 4  & 30 & 5.6 & 1.7   & 1.9    & 0.1 \\
U-R & 4  & 30 & 2.9 & 3.6   & 27.8   & 0.1 \\
U-R & 4  & 30 & 0.7 & 0.4   & 1.0    & 0.0 \\
U-R & 8  & 20 & 6.2 & 3.2   & 3.0    & 0.2 \\
U-R & 8  & 20 & 3.4 & 7.3   & 14.5   & 0.1 \\
U-R & 8  & 20 & 0.9 & 2.8   & 5.7    & 0.0 \\
U-R & 8  & 30 & 6.2 & 4.3   & 4.9    & 0.3 \\
U-R & 8  & 30 & 3.4 & 8.4   & 31.3   & 0.1 \\
U-R & 8  & 30 & 0.9 & 2.2   & 3.9    & 0.0 \\
U-R & 16 & 20 & 6.8 & 6.0   & 5.9    & 0.4 \\
U-R & 16 & 20 & 3.6 & 16.4  & 22     & 0.2 \\
U-R & 16 & 20 & 0.9 & 14.8  & 26.5   & 0.1 \\
U-R & 16 & 30 & 6.8 & 10.0  & 11.3   & 0.4 \\
U-R & 16 & 30 & 3.6 & 24.3  & 101.9  & 0.3 \\
U-R & 16 & 30 & 0.9 & 7.6   & 13.7   & 0.1 \\
U-R & 64 & 20 & 8.0 & 106.2 & 34.5   & 2.3 \\
U-R & 64 & 20 & 4.2 & 129.8 & 221.1  & 0.6 \\
U-R & 64 & 20 & 1.2 & 218   & 286.8  & 0.5 \\
U-R & 64 & 30 & 8.0 & 163.7 & 134.0  & 2.6 \\
U-R & 64 & 30 & 4.2 & 328.6 & 1757.1 & 0.9 \\
U-R & 64 & 30 & 1.2 & 253.3 & 246.3  & 0.5 \\
UN-R & 4  & 20 & 5.6 & 1.1   & 1.1    & 0.1 \\
UN-R & 4  & 20 & 2.9 & 2.4   & 5.8    & 0.0 \\
UN-R & 4  & 20 & 0.7 & 0.5   & 0.9    & 0.0 \\
UN-R & 4  & 30 & 5.6 & 1.8   & 1.8    & 0.1 \\
UN-R & 4  & 30 & 2.9 & 3.2   & 22.7   & 0.1 \\
UN-R & 4  & 30 & 0.7 & 0.4   & 0.9    & 0.0 \\
UN-R & 8  & 20 & 6.2 & 3.7   & 2.6    & 0.2 \\
UN-R & 8  & 20 & 3.4 & 5.8   & 11.6   & 0.1 \\
UN-R & 8  & 20 & 0.9 & 2.6   & 5.8    & 0.0 \\
UN-R & 8  & 30 & 6.2 & 3.8   & 4.0    & 0.3 \\
UN-R & 8  & 30 & 3.4 & 7.2   & 29     & 0.1 \\
UN-R & 8  & 30 & 0.9 & 2.4   & 4.8    & 0.0 \\
UN-R & 16 & 20 & 6.8 & 6.7   & 6.6    & 0.4 \\
UN-R & 16 & 20 & 3.6 & 20.0  & 33.6   & 0.2 \\
UN-R & 16 & 20 & 0.9 & 10.6  & 19.9   & 0.1 \\
UN-R & 16 & 30 & 6.8 & 8.0   & 9.6    & 0.4 \\
UN-R & 16 & 30 & 3.6 & 18.8  & 81.8   & 0.3 \\
UN-R & 16 & 30 & 0.9 & 8.5   & 14.5   & 0.1 \\
UN-R & 64 & 20 & 8   & 108.1 & 40.5   & 2.4 \\
UN-R & 64 & 20 & 4.2 & 239.4 & 468.5  & 0.7 \\
UN-R & 64 & 20 & 1.2 & 212.8 & 358.1  & 0.5 \\
UN-R & 64 & 30 & 8   & 159.4 & 110.6  & 2.5 \\
UN-R & 64 & 30 & 4.2 & 254.8 & 946.1  & 0.8 \\
UN-R & 64 & 30 & 1.2 & 314.7 & 402.1  & 0.4\\
\hline
\end{tabular}
\caption{Computing time CAC  vs Standard}
\label{tab:com_baseline}
\end{table}

\end{document}